\documentclass{amsart}

\usepackage[english]{babel}
\usepackage{amsmath,amssymb}
\usepackage[toc,page]{appendix}

\numberwithin{equation}{section}

\newtheorem{thm}{Theorem}[section]
\newtheorem{prop}[thm]{Proposition}
\newtheorem{lem}[thm]{Lemma}
\newtheorem{cor}[thm]{Corollary}

\newtheorem{oss}{Remark}
\newtheorem*{step}{Step}
\newtheorem{defn}{Definition}[section]
\def\tr{\text{tr}}
\def\eps{\varepsilon}
\def\d{\delta}
\newcommand{\Z}{{\mathbb Z}}
\newcommand{\Si}{\mathbf{S}}
\newcommand{\R}{{\mathbb R}}
\newcommand{\Mi}{{\mathbf M}}
\begin{document}

\title[Large deviations for fast volatility models]{Large deviations  for some fast
 stochastic volatility models by viscosity methods  }
 \author{Martino Bardi , Annalisa Cesaroni , Daria Ghilli   }
 \address{Dipartimento di Matematica,  Universit\`a
di Padova, via Trieste 63, 35121 Padova, Italy}
\email{bardi@math.unipd.it,   \  acesar@math.unipd.it \  ghilli@math.unipd.it}
 \thanks{Partially supported by 
the Fondazione CaRiPaRo Project "Nonlinear Partial Differential Equations: models, analysis, and control-theoretic problems" 
 and the European Project Marie Curie ITN "SADCO - Sensitivity Analysis for Deterministic Controller Design". }
 
\begin{abstract}
We consider the short time behaviour of stochastic systems affected by a stochastic volatility evolving at a faster time scale. We study the asymptotics of a logarithmic functional of the process by methods of the theory of homogenisation and singular perturbations for fully nonlinear PDEs. We point out three regimes depending on how fast the volatility oscillates relative to the horizon length. We prove a large deviation principle for each regime and apply it to 
the asymptotics of option prices near maturity.
\end{abstract}

\maketitle 
\section{Introduction}
In this paper we are interested in 
 stochastic differential equations 
  with two  small parameters $\eps >0$ and $\delta >0$ of the form
\begin{equation}\label{sistema0}
\left\{
\begin{array}{ll}
d X_t =   \eps  \phi(X_t, Y_t) dt + \sqrt{2\eps }\sigma(X_t, Y_t) dW_t \quad &X_{
0
}=x
 \in \R^n,\\
dY_t=\frac{\eps }{\delta} b(Y_t) dt + \sqrt{\frac{2\eps}{\delta}}\tau(Y_t) dW_t \quad &Y_{
0
}=y
 \in \R^m,
\end{array}
\right.\,
\end{equation}
where 
$W_t$ is a standard  $r$-dimensional Brownian motion, the functions $\phi(x,y)$, $\sigma(x,y), b(y),\tau(y)$ are $\Z^m$-periodic with respect to the 
 variable $y$, and the matrix $\tau$ is non-degenerate. 
This  is a model of systems where 
the variables $Y_t$ 
  evolve at a much faster time scale $
 s=\frac{t}{\delta}$ than the other variables $X_t$. The second parameter $\eps$ is added in order to study the  small time behavior of the system, in particular the time has been rescaled in \eqref{sistema0} as $ t \mapsto \eps t$.
Passing to the limit as $\delta\to 0$, with $\eps$ fixed, is a classical singular perturbation problem, its solution leads to the elimination of the state variable $Y_t$ and to the definition of an averaged system defined in $\R^n$ only. There is a large literature on the subject, see the monographs  \cite{K1}, \cite{kp},  the memoir \cite{BA2} and the references therein.  
Here we study the asymptotics as 
 both 
  parameters go to $0$ 
  and we expect different limit behaviors depending on the rate $\eps/\delta$. 
 Therefore we put 
 $$\delta=\eps^\alpha ,
  \; \text{with }\, 
  \alpha>1 ,
  $$ and consider a functional of the trajectories of \eqref{sistema0} of the 
  form
\begin{equation}\label{eqn:veps}
v^\eps(t,x,y) := \eps \log E\left[e^{
h(X_{t})/\eps} | 
 (X.,Y.)\, \, \mbox{satisfy \eqref{sistema0}}\right],
\end{equation}
where $ h \in BC(\mathbb{R}^n)$. The logarithmic form of this payoff is motivated by the applications to large deviations that we want to give. It is known that $v^\eps$ solves the Cauchy problem with initial data  $v^\eps(0,x,y)=h(x)$ for a fully nonlinear parabolic equation. Letting $\eps\to 0$ in this PDE is a regular perturbation of a singular perturbation problem, for which we can rely on the techniques of \cite{ABM}, stemming from Evans' perturbed test function method for homogenisation \cite{E1} and its extensions to singular perturbations \cite{ab, BA, BA2}. 
 We show that under suitable assumptions the functions $v^\eps(t,x,y)$ converge 
  to a function $v(t,x)$
characterised as the solution of the Cauchy problem for a first order Hamilton-Jacobi equation 
\begin{equation}\label{eqn:limit}
v_t-\bar H(x,Dv)=0 \;\text{ in } ]0,T[\times \R^n,\quad v(0,x)=h(x) .
\end{equation}
A significant part of the paper is devoted to the analysis of the \emph{effective Hamiltonian} $\bar H$, which is obtained by solving a suitable cell problem. 
As usual in the theory of homogenisation for fully nonlinear PDEs, this is an additive eigenvalue problem. It turns out  to have different forms in the following three 
 regimes depending on $\alpha$:
$$
\left\{
\begin{array}{lll}
\alpha > 2 \quad \mbox{supercritical case,}\\ 
\alpha=2 \quad \mbox{critical case,}\\ 
\alpha<2 \quad \mbox{subcritical case}. 
\end{array}
\right.\,
$$
More precisely, in the supercritical case the cell problem involves a linear elliptic operator and $\bar H$ has the explicit formula 
\[
\bar{H}(
{x},
{p}) = \int_{\mathbb{T}^m} \! |\sigma(
{x}, y) ^T 
{p}|^2 \, d\mu(y) 
\]
where $\mu$ is the invariant probability measure on the $m$-dimensional torus $\mathbb{T}^m$ of the stochastic process
 \[
 dY_t=b(Y_t)dt+\sqrt{2}\tau(Y_t)dW_t.
\]
In the critical case the cell problem is a fully nonlinear elliptic PDE and $\bar H$ can be represented in various ways based, e.g., on stochastic control.
Finally, in the subcritical case the cell problem is of first order and nonlinear, and a representation formula for   $\bar H$ can be given in terms of deterministic control. In particular, under the condition $\tau\sigma^T=0$ of non-correlations among the components of the white noise acting on the slow and the fast variables in \eqref{sistema0}, we have 
\[
\bar{H}(
{x}, 
{p})= \max_{y\in\R^m}|\sigma^T(
{x},y)
{p}|^2 .
\]
Let us mention that an important step of the method is the comparison principle for the limit Cauchy problem \eqref{eqn:limit},  ensuring that the weak convergence of the relaxed semilimits is indeed uniform, as well as the uniqueness of the limit. It is known that this property of the effective Hamiltonian may require additional conditions \cite{BA2}. Here we show that no extra assumptions are needed in the super- and subcritical cases, whereas in the critical case
the comparison principle holds if either the matrix $\sigma$ is independent on $x$, or it is non-degenerate, or the non-correlation condition $\tau\sigma^T=0$  holds.

The main application of the convergence results is a large deviations analysis of \eqref{sistema0} in the three different regimes. We prove that the measures associated to the process $X_t$ in \eqref{sistema0} satisfy a Large Deviation Principle (briefly, LDP) with good rate function
\begin{equation*}
I(x;x_0,t):= \inf\left[\int_0^t \! \bar{L}\left(\xi(s),\dot{\xi}(s)\right)\, ds \ \Big|\ \xi\in AC(0,t), \ \xi(0)=x_0, \xi(t)=x\right],
\end{equation*}
where $\bar{L}$ is the \emph{effective Lagrangian} associated to $\bar H$ via convex duality. In particular we get that
\begin{equation*}
P(X^\eps_t \in B)=e^
{-\inf_{x \in B} \frac{I(x;x_0,t)}{\eps} + o(\frac 1{\eps})},\; \text{ as } \eps\to 0 
\end{equation*}
for any open set $B\subseteq\R^n$.
Following \cite{FFK} we also apply 
this result to 
an 
 estimate of option prices  near maturity and an asymptotic formula for the implied volatility. 

\smallskip
Our first motivation for the study of systems of the form \eqref{sistema0} comes from financial models with stochastic volatility. 
In such models the vector $X_t$ represents the log-prices of $n$ assets (under a risk-neutral 
probability measure) 
whose 
 volatility $\sigma$ 
 is affected by a process 
 $Y_t$ driven by another Brownian motion, which is often negatively correlated  with the one driving the stock prices (this is the empirically observed leverage effect, i.e., asset prices tend to go down as volatility goes up). Fouque, Papanicolaou, and  Sircar 
 argued in   \cite{FPS} that the bursty  behaviour of volatility observed in financial markets can be described by introducing a faster time scale for a mean-reverting process $Y_t$ by means of the small parameter $\delta$ in \eqref{sistema0}. Several extensions, applications to a variety of financial problems, and rigorous justifications of the asymptotics can be found in 
 \cite{fps1, fps2, BCM, bc, fpss}, 
  see also  the references therein.  On the other hand, Avellaneda et al. \cite{ABBF2} used the theory of large deviations to give asymptotic estimates for the Black-Scholes implied volatility of option prices near maturity in models with constant volatility. 
  In the recent paper \cite{FFK}, Feng, Fouque, and Kumar study the large deviations for system of the form 
 \eqref{sistema0} in the one-dimensional case  $n=m=1$, assuming that $Y_t$ is an Ornstein-Uhlenbeck process and the coefficients in the equation for $X_t$ do not depend on $X_t$. In their model $\eps$ represents a short maturity for the options, $1/\delta$ is the rate of mean reversion of $Y_t$, and the 
 asymptotic analysis is performed for $\delta=\eps^\alpha$ in the regimes $\alpha=2$ and $\alpha=4$. Their methods are based on the approach to large deviations developed in \cite{FK}. A related paper is \cite{FFF} where the Heston model was studied in the regime $\delta=\eps^2$ by methods different from \cite{FFK}.

Although sharing some motivations with \cite{FFK} our results are quite different: 
we treat vector-valued processes under rather general conditions and discuss all the regimes depending on the parameter $\alpha$; our methods are also different, mostly from the theory of viscosity solutions for fully nonlinear PDEs 
and from the theory of homogenisation and singular perturbations for such equations. Our assumption of periodicity with respect to  the $y$ variables may sound restrictive for the financial applications. It is made mostly for technical simplicity and can be relaxed to the ergodicity of the process $Y_t$ as in \cite{BCM, bc}: this will be treated in a paper in preparation.

Large deviation principles have a large literature for diffusions with vanishing noise; some of them were extended 
to two-scale systems with small noise in the slow variables,  
see  
\cite{LIP},  \cite{v}, and more recently  \cite{k2}, \cite{DS}, and \cite{SPI}.  Our methods can be also applied to this different scaling. The paper by Spiliopoulos \cite{SPI} also states some results for the scaling of \eqref{sistema0} under the  assumptions of periodicity and $n=m=1$, but its methods based on weak convergence are completely different from ours. A related paper on 
 homogenisation of a fully nonlinear PDE with vanishing viscosity is \cite{
 CCM}.

The paper is organized as follows. In Section \ref{2} we give the precise assumptions and describe the parabolic PDEs satisfied by $v^\eps$ in the different regimes. 
In sections \ref{3},  \ref{sec:casosupercritico}, \ref{5} 
we analyse the cell problem and the properties of the effective Hamiltonian in the critical 
($\alpha=2$),  supercritical 
 ( $\alpha>2$), and  subcritical case ($\alpha<2$), respectively. 
  Section \ref{6} 
 is devoted to the convergence result for each regime of the functions \eqref{eqn:veps} to the unique viscosity solution 
 of the limit problem \eqref{eqn:limit} with $\bar H$ identified in the previous sections, see Theorems \ref{thm:teoconv1} and 
  \ref{thm:teoconv2}.  
 In section \ref{7} we prove 
the  Large Deviation Principle for all the regimes, Theorem \ref{thm:ldp}. 
 Finally, in Section \ref{8} we give some applications to option pricing.

\section{The fast  stochastic volatility problem}\label{2}
\subsection{The stochastic volatility model}\label{asssistema}
We consider fast-mean reverting stochastic volatility system that can be written in the form
\vspace{0.5cm}
\begin{equation}\label{sistema1}
\left\{
\begin{array}{ll}
d X_t = \phi(X_t , Y_t ) dt + \sqrt{2}\sigma(X_t , Y_t )  dW_t, \quad & X_0 =x \in \R^n\\
d Y_t = \eps^{-\alpha}b(Y_t ) dt + \sqrt{ 2 \eps^{-\alpha}} \tau(Y_t ) dW_t,  \quad & Y_0 =y \in \R^m.
\end{array}
\right.\,
\end{equation}
where $\eps >0$, $ \alpha>1$ and $W_t$ is an $r$-dimensional standard Brownian motion.  We assume  $\phi:\R^n \times \R^m \rightarrow \mathbb{R}^n, \sigma:\mathbb{R}^n \times \mathbb{R}^m \rightarrow \Mi^{n,r}$ are bounded continuous functions, Lipschitz continuous in $(x,y)$ and periodic in $y$, where
$\Mi^{n,r}$ denotes the  set of $n\times r$ matrices. Moreover  $b:\R^m\to\R^m, \tau:\R^m\to \Mi^{m,r}$ are locally Lipschitz continuous functions, periodic in $y$. These assumptions will hold throughout the paper. 
We will use the symbol $\Si^k$ to denote the set of $k\times k$ symmetric matrices. \\
In the following we will assume the uniform nondegeneracy of the diffusion driving the fast variable $Y_t$, i.e for some $\theta>0$
\begin{equation} \label{tau} 
 \xi ^T\tau(y) \tau(y)^T 
  \xi =|\tau^T(y) \xi |^2 >\theta|\xi|^2\quad \mbox{for every}\, \, y\in \mathbb{R} , \xi \in \mathbb{R}^m.
\end{equation}
In order to study small time behavior of the system \eqref{sistema1}, we rescale time $t \rightarrow \eps t$ for $0<\eps \ll 1$, so that the typical maturity will be of order of $\epsilon$. Denoting the rescaled  processes by $X^\eps_t$ and $Y^\eps_t$ we get 
\vspace{0.5cm}
\begin{equation}\label{eqn:systemscaled}\begin{cases} 
d X^\eps_t = \eps\phi(X^\eps_t , Y^\eps_t ) dt + \sqrt{2\eps}\sigma(X^\eps_t , Y^\eps_t )  dW_t,  & 
X^\eps_0 =x \in \R^n\\
d Y^\eps_t =\eps^{1-\alpha }b(Y^\eps_t ) dt + \sqrt{2\eps^{1-\alpha }} \tau(Y^\eps_t ) dW_t,    & Y^\eps_0 =y \in \R^m.
\end{cases}
\vspace{0.5cm}
\end{equation}
 
Next we consider the   functional

\begin{equation}\label{eqn:ueps}
u^{\eps}(t,x,y) := E\left[g(X_{t}) \,|\, 
 (X^\eps.,Y.^\eps)\, \, \mbox{satisfy \eqref{eqn:systemscaled}}\right]  
\end{equation}
where $g \in BC(\mathbb{R}^n)$. We denote with $BC(\mathbb{R}^n)$ the space of bounded continuous functions in $\R^n$. 

The partial differential equation associated to the functions $u^\eps$ 
is 
\begin{multline}
\label{eqn:equazione}
u_t -\eps\tr(\sigma\sigma^{T}D^2_{xx}u) - \eps \phi \cdot D_xu - 2 \eps^{1-\frac{\alpha}{2}}\tr(\sigma\tau^T D^2_{xy}u
)\\ - \eps^{1-\alpha}b\cdot D_yu  -\eps^{1-\alpha}\tr(\tau \tau^T D^2_{yy}u) =0
\end{multline}
in $(0,T) \times \mathbb{R}^n \times \mathbb{R}^m$, where $b$ and $\tau$ are computed in $y$, $\phi$ and $\sigma$ are computed in $(x,y)$. 
The equation is complemented with the initial condition:
$$
u(0,x,y)=g(x).
$$

\begin{oss}\label{oss:per}\upshape
Note that, since we assume the periodicity in $y$ of the coefficients of the equation $b, \sigma, \tau, \phi,$ we have that the solution $u^{\eps}$ of the equation \eqref{eqn:equazione} is periodic in $y$ itself.
\end{oss}

\subsection{The log-tranform and its HJB equation}\label{log-tran}
We introduce the logarithmic transformation method (see \cite{FS}). Assume that 
$$g(x)=e^{
h(x)/\eps} \;\text{ with } h \in BC(\mathbb{R}^n)
$$ 
and define 
\begin{equation}\label{v-eps}
v^{\eps}(t,x,y): = \eps \log u^{\eps}= \eps \log E\left[e^{
h(X^\eps_{t})/\eps}\, | \, 
(X^\eps.,Y^\eps.)\, \, \mbox{satisfy \eqref{eqn:systemscaled}
}\right],
\end{equation}
where $u^\eps$ is defined in \eqref{eqn:ueps}, $x\in\R^n$, $y\in\R^m$, and $t\geq 0$. By \eqref{eqn:equazione} and some computations 
one sees that the equation associated to $v^\eps$ is 
\begin{multline}
\label{eqn:equazione3}
 v_{t} 
 =  |\sigma^T D_xv|^2 + \eps \tr (\sigma\sigma^T D^2_{xx}v) + \eps\phi \cdot D_xv + 2\eps^{-\frac{\alpha}{2}} (\tau \sigma^ T D_xv ) \cdot D_yv+\\ 
 2\eps^{1-\frac{\alpha}{2}}\tr(\sigma\tau^T D^2_{xy}v)  + \eps^{1-\alpha}b\cdot D_y v + \eps^{-\alpha}|\tau^T D_y v|^2 + \eps^{1-\alpha}\tr(\tau \tau^T D_{yy}^2 v), 
\end{multline}
 where $b$ and $\tau$ are computed in $y$, $\phi$ and $\sigma$ are computed in $(x,y)$. 
In general, the functions $u^\epsilon$ are not smooth but 
one can check  that  $v^\epsilon$ is a viscosity solutions of \eqref{eqn:equazione3} 
 (see in particular Chapter VI and VII of \cite{FS}). 

In the following proposition we characterize the value function $v^\eps$ as the unique continuous viscosity solution to a suitable parabolic problem with initial data for each of the three regimes. 
A general reference for these issue is \cite{FS}. 
The equation \eqref{eqn:equazione3} satisfied by $v^\epsilon$ 
involves a quadratic nonlinearity in the gradient. 
This case was studied by Da Lio and Ley in \cite{DLL}, where the reader can find a proof of the next result.
\begin{prop}\label{prop:propequazione} 
i) Let  $\alpha \geq2$ and define 
\begin{eqnarray*} 
H^{\eps}(x,y, p, q, X,Y, Z)&:=&   |\sigma^T p|^2+ b\cdot q+ \tr(\tau \tau^T Y) +  \eps\left(\tr(\sigma\sigma^T X) + \phi \cdot p\right) 
\\ &+& 2\eps^{\frac{\alpha}{2}-1} (\tau \sigma^T p) \cdot q   +2\eps^{\frac{1}{2}}\tr(\sigma\tau^T Z) + \eps^{\alpha-2}|\tau^T q|^2  .
\end{eqnarray*}
Then $v^\eps$  
 is the unique bounded continuous viscosity solution of the 
  Cauchy problem
\begin{equation}\label{eqn:equazioneeps1}
\begin{cases} 
 \partial_t v^\eps -H^\eps \left(x,y,D_x v^\eps, \frac{D_y v^\eps}{\eps^{\alpha-1}}, D^2_{xx} v^\eps, \frac{D^2_{yy} v^\eps}{\eps^{\alpha-1}}, \frac{D^2_{xy} v^\eps}{\eps^{\frac{\alpha-1}{2}}}\right)=0 & \, \, \mbox{in} \, \, [0,T] \times \mathbb{R}^n \times \mathbb{R}^m,\\
v^\eps(0,x,y)=h(x) & \mbox{ in }  \mathbb{R}^n\times \mathbb{R}^m. 
\end{cases}
\end{equation}

ii) Let $\alpha <2$ and define 
\begin{eqnarray*}
H_{\eps}(x,y, p, q, X,Y, Z)&:=& |\sigma^T p|^2+  |\tau^T q|^2+2(\tau \sigma^T  p)\cdot q +  \eps\left(\tr(\sigma\sigma^T X) + \phi\cdot p\right) \\ & +&  \eps^{1-\frac{\alpha}{2}}(b\cdot q +\tr(\tau \tau^T Y))  +2\eps^{1-\frac{\alpha}{4}}\tr(\sigma\tau^T Z).\end{eqnarray*}
Then $v^\eps$  
 is the unique bounded continous viscosity solution of the 
  Cauchy problem
\begin{equation}\label{eqn:equazioneeps2}
\begin{cases} 
 \partial_t v^\eps -H_\eps \left(x,y,D_x v^\eps, \frac{D_y v^\eps}{\eps^{\frac{\alpha}{2}}}, D^2_{xx} v^\eps, \frac{D^2_{yy} v^\eps}{\eps^{\frac{\alpha}{2}}}, \frac{D^2_{xy} v^\eps}{\eps^{\frac{\alpha}{4}}}\right)=0  &  \mbox{ in  }  [0,T] \times \mathbb{R}^n \times \mathbb{R}^m,\\
v^\eps(0,x,y)=h(x) & \, \, \mbox{in} \, \, \mathbb{R}^n\times \mathbb{R}^m. 
\end{cases}
\end{equation}
\end{prop}

Our goal is to study the limit as $\eps\to 0$ of the functions $v^\eps$ described in Proposition \ref{prop:propequazione}. Following the viscosity solution apoproach to singular perturbation problems (see  \cite{BA2},\cite{BA}), we  define a limit or effective Hamiltonian $\overline{H}$ and we characterize the limit of $v^\eps$ as the unique solution of an appropriate Cauchy problem with Hamiltonian $\overline{H}$. The first step in the procedure is the identification of the limit Hamiltonian. In order to define this operator, we make the 
 ansatz that 
 the function  $v^{\eps}$ admits 
 the formal asymptotic expansion 
\begin{equation}\label{eqn:exp}
v^{\eps}(t,x,y)=v^{0}(t,x) + \eps^{\alpha-1} w(t,x,y)
\end{equation}
and plug it 
 into the equation. 
In the following sections we show that the limit Hamiltonian is different in the  three different regimes: the critical case ($\alpha=2$), the supercritical case (when $\alpha >2$), and the subcritical case (when $\alpha <2$). 

  Numerical experiments in \cite{TY}   indicate that the first order approximation in the 
   expansion \eqref{eqn:exp} is sufficiently accurate 
    to find option
prices in a fast mean-reversion case of the volatility process.  

\section{The critical case: $\alpha=2$}
\label{sec:casocritico}
\label{3}

Equation \eqref{eqn:equazione3} with $\alpha=2$ becomes 
\begin
{align}
\label{eqn:equazione3r2}
v_{t} = |\sigma^T D_xv|^2+ \eps \tr(\sigma\sigma^T D^2_{xx}v) + \eps\phi \cdot D_xv + \frac{2}{\eps}  (\tau \sigma^ T D_x v ) \cdot D_y v\\-2\tr(\sigma\tau^T D^2_{xy}v)
+\frac{1}{\eps}b\cdot D_y v +\frac{1}{\eps^2}| \tau^T D_y v|^2 + \frac{1}{\eps}\tr(\tau \tau^T D_{yy}^2 v).\notag
\end
{align}
\subsection{The effective Hamiltonian}
We plug in the equation \eqref{eqn:equazione3r2} the formal asymptotic expansion 
$$
v^{\eps}(t,x,y)=v^{0}(t,x) + \eps  w(t,x,y)
$$ 
and we  obtain  
$$
v^0_{t} -|\sigma^T D_xv^0|^2 - 2 (\tau \sigma^ T D_x v^0 ) \cdot D_y w  - b\cdot D_y w -|\tau^T D_y w|^2 - \tr(\tau \tau^T D_{yy}^2 w) =O(\eps).
$$
We want to eliminate the \textit{corrector} $w$ and the dependence on $y$ in this equation and remain with a left 
hand side of the form  $v^0_{t}-\bar{H}({x}, D_x v^0)$. 
 Therefore we freeze 
$\bar{x}$ and $\bar{p}=D_xv^0(\bar x)$ and define the \textit{effective Hamiltonian} $\bar{H}(\bar{x}, \bar{p})$ as  the unique
constant  such that the following stationary PDE in $\R^m$, called \textit{cell problem},  has a viscosity solution $w$: 
\begin{equation}\label{eqn:cellacritico}
\bar{H}(\bar{x},\bar{p})- |\sigma^T\bar{p}|^2 - 2 (\tau \sigma^ T \bar{p} ) \cdot D_y w(y)  - b\cdot D_y w(y) -|\tau^T D_y w(y)|^2 - \tr(\tau \tau^T D_{yy}^2 w(y))=0,
\end{equation}
where $\sigma$ is computed in $(\bar{x},y)$ and $\tau, b$ in $y$. This is an additive eigenvalue problem that arises the theory of ergodic control and has a wide literature.
Under our standing assumptions we have  the following result.
\begin{prop}\label{thm:trucell}
For any fixed $(\bar{x}, \bar{p})$, there exists a unique   $\bar{H}(\bar{x}, \bar{p})$ for which the equation \eqref{eqn:cellacritico} has a  periodic viscosity solution $w$. Moreover 
 $w \in C^{2,\alpha}$ for some $0<\alpha<1$ and satisfies for some $C>0$ independent of $\bar{p}$ and $\forall \bar{x},\bar{p} \in \mathbb{R}^n$
\begin{equation}\label{eqn:stimagrad}
\max_{y\in \mathbb{R}^m}|Dw(y;\bar{x},\bar{p})|\leq C(1+|\bar{p}|).
\end{equation}
\end{prop}
To prove Proposition \ref{thm:trucell}, we need the following lemma.
\begin{lem}\label{stimagrad}
Let $\delta >0$ and $w_\delta(\cdot;\bar{x},\bar{p}) \in C^2(\mathbb{R}^m)$ be a periodic solution of 
\begin{equation}\label{cell}
\delta w_\delta +F(\bar{x},y, \bar{p},Dw_\delta, D^2w_\delta)-|\sigma(\bar{x},y)\bar{p}|^2=0,
\end{equation}
where 
\begin{equation}\label{F}
 F(\bar x,y,\bar p, q,Y):= - \tr(\tau \tau^T(y) Y)-|\tau^T(y) q|^2 -b(y)\cdot q -2 (\tau(y) \sigma^T(\bar x,y) \bar{p} ) \cdot q .
  \end{equation}
 Then there exists $C>0$ independent of $\bar{p}$  such that for all $\bar{x}, \bar{p} \in \mathbb{R}^n$ it holds
\begin{equation}\label{bern}
\max_{y\in \mathbb{R}^m}|D_y w_\delta(y;\bar{x},\bar{p})| \leq C(1+|\bar{p}|).
\end{equation}
\end{lem}
\begin{proof}
The proof uses the Bernstein method, following  the derivation of similar estimates  in \cite{EI}. 
We carry out the computations  in the case $\tau,\sigma,b$ are $C^1$. When $\tau,\sigma, b$ are Lipschitz  the result can be proved by smooth approximation.\\
Denote by $w^\delta:=w_\delta(y;\bar{x},\bar{p})$ the solution of \eqref{cell}. 
By comparison 
 with constant sub- and supersolutions we get the uniform bound
\begin{equation}\label{wk}
|\delta w^\d|\leq \max_{y\in\mathbb{R}^m}|\sigma^T(\bar{x},y)\bar{p}|^2\quad \forall y \in \mathbb{R}^m.
\end{equation}
Define the function $z$ as follows
$$
z:=|D w^\d|^2.
$$
Should $z$ attains its maximum at some point $y_0$, then at $y_0$
\begin{equation}\label{uno}
z_i=2 w^\d_{k} w^\d_{ki}=0 \quad i=1,\dots,m,
\end{equation}
where we are adopting the summation convention, and
\begin{equation}\label{due}
0 \leq -(\tau \tau^T)_{ij} z_{ij}=-2(\tau \tau^T)_{ij} w^\d_{ki}w^\d_{kj} - 2w^\d_k(\tau \tau^T)_{ij}w^\d_{ijk}.
\end{equation}
Then at $y_0$
\begin{multline*}
\theta |D^2w^\d|^2 \leq (\tau \tau^T)_{ij} w^\d_{ki} w^\d_{kj} \leq\\ - w^\d_k(\tau\tau^T)_{ij} w^\d_{ijk}=-w^\d_k\left((\tau\tau^T)_{ij} w^\d_{ij}\right)_k + w^\d_k(\tau\tau^T)_{ij,k} w^\d_{ij},
\end{multline*}
where we have used \eqref{due}.
Thus at $y_0$
\begin{multline*}
\theta |D^2w^\d|^2 \leq\\ w^\d_k\left(-\delta w^\d +(2\tau\sigma^T\bar{p}+b) \cdot Dw^\d +|\tau^T D w^\d|^2 +|\sigma^T \bar{p}|^2\right)_k+w^\d_k(\tau\tau^T)_{ij,k} w^\d_{ij},
\end{multline*}
where we have used \eqref{cell}. 
Thanks to \eqref{uno}
\begin{multline*}
w_k^\d(|\tau^TDw^\d|^2)_k=w^\d_k((\tau\tau^T)_{ij}w^\d_iw^\d_j)_k=\\ w^\d_k(\tau\tau^T)_{ij,k}w^\d_iw^\d_j +w_k^\d(\tau \tau^T)_{ij}w^\d_{ik}w^\d_j +w_k^\d(\tau\tau^T)_{ij} w^\d_i w^\d_{jk}=w_k^\d(\tau\tau^T)_{ij,k}w^\d_iw^\d_j.
\end{multline*}
Moreover \[w^\d_k(\tau\tau^T)_{ij,k} w^\d_{ij}\leq \frac{\theta}{2} |D^2w^\d|^2+ \frac{C}{2\theta}  |Dw^\d|^2.\] 
Then 
$$
\theta|D^2w^\d|^2 \leq C(1 +|\bar{p}|)|Dw^\d|^2+C|Dw^\d|^3 +\frac{\theta}{2}|D^2w^\d|^2+C|\bar{p}|^2|Dw^\d|\quad \mbox{at} \,\, y_0
$$
and $C>0$ depends only on the $L^\infty$ norm of $\sigma, b, \tau$ and on the derivatives of $\sigma, b$ and $\tau$. Therefore
\begin{equation}\label{d22}
|D^2w^\d|^2 \leq C(1 +|Dw^\d|^2+|\bar{p}||Dw^\d|^2+|\bar{p}|^2|Dw^\d|^2+|Dw^\d|^3)\quad \mbox{at} \,\, y_0.
\end{equation}
Thanks to the uniform ellipticity of $\tau$ and using equation \eqref{cell}, we have
\[
\begin{split}
\theta |Dw^\d|^2 \leq |\tau^T Dw^\d|^2 = \delta w^\d-\tr(\tau \tau^T D^2 w^\d) -2\tau\sigma^T \bar{p} \cdot D w^\d -b\cdot Dw^\d\quad \mbox{at} \,\, y_0.
\end{split}
\]

Using \eqref{wk}, we get at $y_0$
\begin{equation}\label{z2}
\begin{split}
z^2=|Dw^\d|^4&\leq C(|\bar{p}|^4+|D^2w^\d|^2+|\bar{p}|^2|Dw^\d|^2+|Dw^\d|^2+|\bar{p}||Dw^\d|^2+|\bar{p}|^2|D^2w^\d|\\&+|\bar{p}|^2|Dw^\d|+|\bar{p}|^3|Dw^\d|+|D^2w^\d||\bar{p}||Dw^\d|+|D^2w^\d||Dw^\d|).
\end{split}
\end{equation}

Then \eqref{bern} follows by dividing \eqref{z2} by $|Dw^\d|^3$ and noticing that the right member in \eqref{z2} is polynomial of degree $4$ in $|\bar{p}|$ and $|Dw^\d|$. 

\end{proof}

\begin{proof}
We use the methods of \cite{AL}
 based on the small discount approximation
 \begin{equation}\label{discount}
 \d w_\d +F(\bar x,y,\bar p, D_yw_\d,D^2_{yy}w_\d)-|\sigma^T(\bar x,y)\bar{p}|^2=0 \quad\text{in }\R^m ,
 \end{equation}
 where $F$ is defined in \eqref{F}. 
Let $w_\d:=w_\d(y,\bar{x},\bar{p}) \in C^2(\mathbb{R}^m)$ be a solution of \eqref{discount}. We show that $\d w_\d(y)$ converges along a subsequence of $\d\to 0$ to the constant $\bar{H}(\bar{x}, \bar{p})$ and $w_\d(y)- w_\d(0)$  converges to the corrector $w$. The hard part is proving equicontinuity estimates for $\d w_\d$. Different from \cite{AL, BA2}, here the leading term in \eqref{eqn:cellacritico} is $|\tau^T Dw|^2$ rather than $\tr(\tau \tau^T D^2 w)$. Then the Krylov-Safonov estimates for elliptic PDEs must be replaced by the Lipschitz estimates proved in Lemma \ref{stimagrad}. In fact, thanks to \eqref{bern}, for some $C>0$ independent of $\bar{p}$ and for all $y,z \in \mathbb{R}^m$ and $\d>0$
\begin{equation}\label{equilip}
| \d w_\delta(y)- \d w_\d(z)|\leq  C\d (1+|\bar{p}|)|y-z| 
\end{equation}
and the equicontinuity follows. The equiboundness follows from \eqref{wk}. Then by Ascoli-Arzela theorem, there is a sequence $\delta_n \rightarrow 0$ such that $\delta_n w_{\delta_n}$ converges locally uniformly to a constant thanks to \eqref{equilip}. We call it $\bar{H}$. 
 Similarly, we prove that $v_\d:=w_\delta(y)-w_\d(0)$ is equibounded and equicontinuous and thus converges locally uniformly along a subsequence   to a function $w$. Then, from \eqref{discount} we get
$$
\d v_\d +\d w_\d(0)+F(\bar x,y,\bar p, D_yv_\d,D^2_{yy}v_\d)-|\sigma^T(\bar x,y)\bar{p}|^2=0 ,\quad\text{in }\R^m .
$$
Since $v_\d $ is equibounded 
 $\delta v_\d \rightarrow 0$. Then from $\d w_\d\rightarrow \bar{H}$ we get that $w$ is a solution of \eqref{eqn:cellacritico}. Finally, by the comparison principle for \eqref{discount}, it is standard to see that $\bar{H}$ is unique.\\
Moreover the regularity theory for viscosity solutions of convex uniformly elliptic equations  
implies that $w\in C^{2,\alpha}$ for some $0<\alpha <1$. \\
Finally the corrector inherits \eqref{equilip} and satisfies for some $C>0$ independent of $\bar{p}$ and for all $\bar{x},\bar{p}\in \mathbb{R}^n$
$$
\max_{y\in \mathbb{R}^m}|D_yw(y;\bar{x},\bar{p})|\leq C(1+|\bar{p}|).
$$

\end{proof}
\subsection{Properties and formulas for $\bar{H}$}
The next result lists  some elementary properties of the effective Hamiltonian $\bar{H}$.
\begin{prop}\label{prop:contbarh1}
\begin{enumerate}\itemsep2pt
\item[(a)] $\bar{H}$  is continuous on $\mathbb{R}^n \times \mathbb{R}^n$;

\item[(b)] the function $p \rightarrow \bar{H}(
{x}, 
{p})$ is convex;

\item[(c)]  
\begin{equation}\label{eqn:bounds}
\min_{y\in\R^m} |\sigma^T(\bar x,y)\bar p|^2\leq \bar H(\bar x,\bar p)\\ \leq \max_{y\in\R^m} |\sigma^T(\bar x,y)\bar p|^2 ;
\end{equation}

\item[(d)] There exists  $C>0$ independent of $p$ such that, for all $x,\bar{x},p \in \mathbb{R}^n$,
\begin{equation}\label{hgrowth}
|\bar{H}(x,p)-\bar{H}(\bar{x},p)|\leq C(1+|p|^2)|x-\bar{x}|;
\end{equation}

 \item[(e)]if
\begin{equation}\label{uncor}
\tau(y)\sigma^T(x,y)=0 \quad \forall x\in \mathbb{R}^n, y \in \mathbb{R}^m,
\end{equation}
then, for all $x, \bar x, p, \bar p\in\R^n$, 
\begin{multline}
\label{est_x_H}
\min_{y\in\R^m}\left(|\sigma^T(x,y)p|^2 - |\sigma^T(\bar x,y)\bar p|^2\right)\leq \bar H(x,p) - \bar H(\bar x,\bar p)\\ \leq \max_{y\in\R^m}\left(|\sigma^T(x,y)p|^2 - |\sigma^T(\bar x,y)\bar p|^2\right).
\end{multline}
\end{enumerate}
\end{prop}

\begin{oss}\rm{
The meaning of assumption \eqref{uncor} is that the components of the Brownian motion $W_t$ influencing the slow variables $X_t$ are not correlated with the components acting on the slow variables $Y_t$. In fact the condition is satisfied if the last $m$ columns of $\sigma$ and the first $n$ columns of $\tau$ are indentically zero.}
\end{oss}

\begin{proof}
The results (a), (b), and (c) are obtained by standard methods in the theory of homogenisation, by means of comparison principles for the approximating equation  \eqref{discount}, see, e.g., \cite{E1, ab}. Let us show one inequality in \eqref{hgrowth} (the other being symmetric).  Let $w_\d(y):=w_\d(y;\bar{x},p)$ ans 
$v_\d(y):=w_\d (y; {x}, {p})$. Then $v_\d$ satisfies

\begin{multline}\label{d}
\d v_\d +F(\bar{x},y,p, D_yv_\d,D^2_{yy}v_\d)-|\sigma^T(\bar x,y)\bar{p}|^2= |\sigma^T(x,y)p|^2 - |\sigma^T(\bar x,y)\bar p|^2+\\(2\tau(y) \sigma^T(x,y) p -2\tau(y) \sigma^T(\bar{x},y)p)\cdot  D v_\d.
\end{multline}
Thanks to Lemma \ref{stimagrad} we estimate $Dv_\d$, and then, using the Lipschitz continuity of $\sigma$, we get for some $C>0$
\begin{equation}
\d v_\d +F(\bar{x},y,p, D_yv_\d,D^2_{yy}v_\d)-|\sigma^T(\bar x,y)\bar{p}|^2\leq C(1+|p|^2)|x-\bar{x}|.
\end{equation}
Then the comparison principle  gives
\[
\d v_\d(y) - \d w_\d(y)\leq C(1+|p|^2)|x-\bar{x}| \quad\forall \,y\in\R^m .
\]
By letting $\d\to 0$ we get  the 
 inequality for $\bar{H}(x,p)-\bar{H}(\bar{x},p)$ in \eqref{hgrowth}, and by exchanging $x$ and $\bar{x}$ we complete the proof.

If \eqref{uncor} holds, \eqref{d} simplifies to
\begin{equation*}
\d v_\d +F(\bar{x},y,p, D_yv_\d,D^2_{yy}v_\d)-|\sigma^T(\bar x,y)\bar{p}|^2\leq \max_y\left\{ |\sigma^T(x,y)p|^2 - |\sigma^T(\bar x,y)\bar p|^2\right\} .
\end{equation*} 
Then, as before, we obtain by comparison the second inequality in \eqref{est_x_H}, and the first is got in a symmetric way.
 
\end{proof}

Next we give some representation formulas for the effective Hamiltonian $\bar{H}$.  
\begin{prop}
(i) $\bar{H}$ satisfies
\begin{equation}\label{eqn:khaise}
\bar{H}(\bar{x}, \bar{p})= \lim_{\delta \rightarrow 0}  \sup_{\beta(\cdot)} \delta E \left[ \int_0^\infty \left (|\sigma(\bar{x}, Z_t)^T \bar{p}|^2 - |\beta(t)|^2\right) e^{-\d t} dt \, |\, Z_0=z\right]
\end{equation}
and 
\begin{equation}\label{khaiset}
\bar{H}(\bar{x},\bar{p})=\lim_{t \rightarrow \infty} \sup_{\beta(\cdot)} \frac{1}{t} E\left[\int_0^{t} (|\sigma^T(\bar{x}, Z_s) \bar{p}|^2 -|\beta(s)|^2) ds\, |\, Z_0=z\right],
\end{equation}
where $\beta(\cdot)$ is an admissible control process 
taking values in $ \mathbb{R}^r$ for the 
stochastic control system 
\begin{equation}\label{Zt}
d Z_t=\left(b(Z_t) +2\tau(Z_t)\sigma^T(\bar{x}, Z_t) \bar{p} -2\tau(Z_t) \beta(t)\right) dt +\sqrt{2} \tau(Z_t) d W_t ; 
\end{equation} 
(ii) moreover  
\begin{equation}\label{meas}
\bar{H}(\bar{x}, \bar{p})=  \int_{\mathbb{T}^m} \left (|\sigma(\bar{x}, z)^T \bar{p}|^2 - |\tau(z)^T Dw(z)|^2\right) d\mu(z),
\end{equation}
where $w=w(\cdot;\bar{x}, \bar{p})$ is the corrector defined in Proposition \ref{thm:trucell} and $\mu=\mu(\cdot;\bar{x}, \bar{p})$ is the invariant probability measure 
of the process  \eqref{Zt1} with the feedback $
{\beta}(z)= -\tau^T (z) D w(z)$;\\
(iii)  finally 
\begin{equation}\label{eqn:barh1}
\bar{H}(\bar{x},\bar{p})=\lim_{t \rightarrow \infty} \frac{1}{t} \log E\left[e^{\int_{0}^{t} \! |\sigma^T(\bar{x}, Y_s) \bar{p}|^2 \, ds \, }|\, Y_{0}=y\right],
\end{equation}
where $Y_t$ is the stochastic process defined by
\begin{equation}\label{eqn:y}
d Y_t=\left(b(Y_t) + 2 \tau(Y_t) \sigma^T(\bar{x},Y_t) \bar{p} \right) dt + \sqrt{2} \tau(Y_t) dW_t.
\end{equation}
\end{prop}

\begin{proof}
(i)  The first formula comes from 
 a control interpretation of the approximating $\delta$-cell problem \eqref{cell}. 
We write it 
as the Hamilton-Jacobi-Bellman equation
\begin{multline}
\label{previous1}
\delta w_\delta  + \\ \inf_{\beta \in \mathbb{R}^r} \left\{ -\tr(\tau(y)\tau(y)^TD^2 w_\d+ \left(2\tau(y) \beta -2\tau(y)\sigma(\bar{x},y)^T \bar{p}-b(y)\right)\cdot D_y w_\delta +|\beta|^2\right\}\\
- |\sigma(\bar{x},y)^T \bar{p}|^2=0
\end{multline}
and we represent $w_\delta$ as the value function of the infinite horizon discounted stochastic control problem  (see, e.g., \cite{FS})
$$
w_\d(z)=\sup_{\beta(\cdot)} E\left[\int_0^{\infty} (|\sigma^T(\bar{x}, Z_t) \bar{p}|^2 -|\beta(t)|^2) e^{-\delta t}dt\, | \, Z_0=z\right],
$$
where $Z_t$ is defined by \eqref{Zt}. Then \eqref{eqn:khaise} follows from the proof of Proposition \ref{thm:trucell}. 
 
For the formula \eqref{khaiset} we consider  the \textit{t-cell problem}
\begin{equation}\label{eqn:cellat}
\left\{
 \begin{array}{ll}
 \frac{\partial v}{\partial t} -\tr(\tau \tau^T D^2 v) - |\tau^T D v|^2 - (b + 2 \tau \sigma^T \bar{p}) \cdot Dv - |\sigma^T \bar{p}|^2=0 &\mbox{in } (0,+\infty)\times \mathbb{R}^m,\\
 v(0,z)=0 &\mbox{on } \mathbb{R}^m.
                 \end{array}
\right.\,
\end{equation}
This is also a HJB equation, whose solution is the value function
\[
v(t,z;\bar{x},\bar{p})=\sup_{\beta(\cdot)}  E\left[\int_0^{t} (|\sigma^T(\bar{x}, Z_s) \bar{p}|^2 -|\beta(s)|^2) ds\, |\, Z_0=z\right],
\]
where $Z_t$ is defined by \eqref{Zt}.
Then a generalized Abelian-Tauberian theorem (see  \cite{BA} for a general proof based only on the comparison principle for the Hamiltonian) states that 
\begin{equation}\label{limcellat}
\bar{H}(\bar{x},\bar{p})=\lim_{t \rightarrow +\infty} \frac{v(t,z;\bar{x},\bar{p})}{t} \quad \mbox{ uniformly in }z.
\end{equation}

(ii) The formula \eqref{meas} is derived from a direct  control interpretation of the cell problem \eqref{eqn:cellacritico}. In fact, it is the HJB equation of the ergodic control problem of maximizing 
\[
\lim_{T \rightarrow \infty} \frac{1}{T} E\left[\int_0^{T} (|\sigma^T(\bar{x}, Z_s) \bar{p}|^2 -|\beta(s)|^2) ds\, |\, Z_0=z\right],
\]
among admissible controls $\beta(\cdot)$ taking values in $\R^r$ for the system \eqref{Zt}, as before.
The process $Z_t$ associated to each control is ergodic with a unique invariant measure $\mu$ on $\mathbb{T}^m$ because it is a nondegenerate diffusion on $\mathbb{T}^m$, see, e.g., \cite{BA2}, so the limit in the payoff functional exists and it is the space average in $d \mu$ of the running payoff. Since the HJB PDE \eqref{eqn:cellacritico} has a smooth solution $w$, it is known from a classical verification theorem that the feedback control that achieves the minimum in the Hamiltonian, i.e.,  ${\beta}(z)= -\tau^T (z) D w(z)$, is optimal. Then \eqref{meas} holds with $\mu$ the invariant measure of the process
\begin{equation}\label{Zt1}
d \tilde{Z}_t=\left(b(\tilde{Z}_t) +2\tau(\tilde{Z}_t)\sigma^T(\bar{x}, \tilde{Z}_t) \bar{p} +2\tau(\tilde{Z}_t) \tau^T(\tilde{Z}_t)Dw(\tilde{Z}_t)\right) dt +\sqrt{2} \tau(\tilde{Z}_t) d W_t.
\end{equation}

(iii) To prove \eqref{eqn:barh1}, take $v=v(t,x;\bar{x},\bar{p})$  a periodic solution of the $t$-cell problem and define the function $f(t,y)= e^{v(t,y)}$. Then $f$ solves the following equation
$$
\left\{
 \begin{array}{ll}
\frac{\partial f}{\partial t} -f|\sigma^T\bar{p}|^2 -(2\tau\sigma^T \bar{p} + b) \cdot D f -\tr(\tau \tau^T D^2 f)=0 &\mbox{in} \, \, (0,\infty) \times \mathbb{R}^m\\
 f(0,z)=1 &\mbox{ in } \mathbb{R}^m.
                 \end{array}
\right.\,
$$
By the Feynman-Kac formula, we have 
$$
f(t,y)=E\left[e^{\int_{0}^t \! |\sigma^T(\bar{x}, Y_s) \bar{p}|^2 \, ds \, }|\, Y_{0}=y\right],
$$
where $Y_t$ is defined by \eqref{eqn:y}. Then
$$
v(t,y)=\log E\left[e^{\int_{0}^t \! |\sigma^T(\bar{x}, Y_s) \bar{p}|^2 \, ds \, }|\, Y_{0}=y\right]
$$
and thanks to  \eqref{limcellat} we get \eqref{eqn:barh1}.
\end{proof}

\begin{oss}\rm{
For $x, p\in\R^n$ define the  following perturbed generator $L^{x,p}$
$$
L^{x,p} g(y) := Lg(y) + 2(\tau \sigma(x,y)^T p) \cdot D_y g(y),
$$
where 
$$
L=b\cdot D_y + \tr(\tau \tau^T D_{yy}^2).
$$
Then the equation \eqref{eqn:cellacritico} becomes
\begin{equation}\label{eqn:bp}
\bar{H}-e^{-w} L^{\bar{x}, \bar{p}} e^{w} - |\sigma^T \bar{p}|^2 = 0,
\end{equation}
because
$
e^{-w} L e^{w}= L w+ |\tau^T D_{y}w|^2
$ 
 gives  
$$
e^{-w} L^{\bar{x},\bar{p}} e^{w}= e^{-w} L e^{w} + 2(\tau \sigma^T \bar{p}) \cdot D_y w= L w +| \tau^T D_{y}w|^2 +2(\tau \sigma^T \bar{p}) \cdot D_y w.
$$
Multiplying \eqref{eqn:bp} by $e^{w}$ we get, for $
g(y)=e^{w(y)} ,
$
\begin{equation}\label{eqn:eigenvalue}
\bar{H}g(y)-(L^{\bar{x},\bar{p}} + V^{\bar{x},\bar{p}})g(y)=0,
\end{equation}
where
$
V^{\bar{x},\bar{p}}(y)=|\sigma^T(\bar{x},y) \bar{p}|^2
$ 
is a multiplicative potential operator. 

We 
conclude  that  if 
 $w$ 
 is a solution of \eqref{eqn:cellacritico}, then $\bar{H}$  \emph{is the first 
  eigenvalue of the linear operator $L^{\bar{x},\bar{p}} + V^{\bar{x},\bar{p}}$, with eigenfunction $g=e^w$.}
}
\end{oss}

\begin{oss} \rm{
Equations like \eqref{eqn:cellacritico} have been studied in an aperiodic setting by Khaise and Sheu in \cite{KSsec}. They prove the existence of a constant $\bar{H}$ such that there is a unique smooth  solution $w$ with prescribed growth of \eqref{eqn:cellacritico}. Moreover they provide a representation formula for $\bar{H}$ as the convex conjugate of  a suitable operator over a space of measures.} 
\end{oss}

\subsection{Comparison principle for $\bar{H}$}
The comparison theorem among viscosity sub- and supersolutions of the 
 limit PDE 
 \begin{equation}\label{eqn:effPDE}
 v_t - \bar{H}(x,Dv)=0 \quad \mbox{in} \, \, (0,T) \times \mathbb{R}^n 
\end{equation}
will be the crucial tool for proving that the convergence of $v^\eps$ is not only in the weak sense of semilimits but in fact uniform, and the limit is unique. It is known from \cite{BA2} that in general the regularity of $\bar{H}$ with respect to $x$ may be worse than that of $H^\eps$ and 
the comparison principle may fail. Next result gives  three alternative additional conditions ensuring the comparison.
\begin{thm}
\label{thm:compbarh1} Assume either one of the following conditions:

\noindent (i) $\sigma$ is independent of $x$, i.e., $\sigma=\sigma(y)$, and $h\in BUC(\mathbb{R}^n)$; or

\noindent (ii) for some $\nu>0$ 
 \begin{equation}\label{coercive}
 |\sigma^T(x,y)p|^2>\nu|p|^2 \quad\forall \, x, p\in\R^n, \,y\in\R^m ;
\end{equation}
or\\
\noindent (iii) $\tau(y)\sigma^T(x,y)=0$ \quad for all $ x \in \mathbb{R}^n,y \in \mathbb{R}^m.$

\noindent  Let $u\in BUSC([0, T]\times\R^n)$ and $v\in BLSC([0,T]\times \R^n)$ be, respectively, a bounded upper semicontinuous subsolution and   a bounded lower semicontinuous supersolution to \eqref{eqn:effPDE}
 such that $u(0,x)\leq h(x)\leq v(0,x)$ for all $x\in \R^n$. Then $u(x,t)\leq v(x,t)$ for all $x\in \R^n$ and $0\leq t\leq T$. 
\end{thm}

\begin{proof}
In case (i) $\bar{H}$ is independent of $x$ and it is continuous by
Proposition \ref{prop:contbarh1}. Then the result follows from standard theory, see e.g. 
\cite{Barl}.  
%

In case (ii) the bounds 
\eqref{eqn:bounds} and \eqref{coercive} give
\[
\bar{H}(
{x}, 
{p})\geq 
\nu |{p}|^2 .
\]
Then $\bar{H}$ is coercive and has the properties (a) and (b) and (c) of Proposition \ref{prop:contbarh1}.  The result follows from \cite{DLL} once we prove that for some $C>0$ and all $x,y,p\in \mathbb{R}^n$
\begin{equation}\label{ddl}
|\bar{L}(x,p)-\bar{L}(y,p)|\leq C(1+|p|^2)|x-y|,
\end{equation}
where $\bar{L}(x,p)$ is the effective Lagrangian, i.e. $\bar{L}(x,p)=\sup_{q\in\mathbb{R}^m}\{p \cdot q-\bar{H}(x,q)\}$. 
 Take $\bar{q}$ such that
$$
\bar{L}(x,p)=\bar{q}\cdot p-\bar{H}(x,\bar{q}).
$$
Then
$$
\bar{L}(x,p)-\bar{L}(y,p)\leq \bar{H}(y,\bar{q}) -\bar{H}(x,\bar{q})\leq C(1+|\bar{q}|^2)|x-y|,
$$
where we have used the property (d) of Proposition \ref{prop:contbarh1}. We want to estimate $|\bar{q}|$. If 
 $|\bar{q}| >\frac{|p|}{\nu}$, then
$$
0\leq \bar{L}(x,p)=\bar{q}\cdot p-\bar{H}(x,\bar{q})\leq \bar{q} \cdot p -\nu|\bar{q}|^2 <0
$$
and we reach a contradiction. Then
$$
\bar{L}(x,p)-\bar{L}(y,p)\leq C\left(1+\frac{|p|^2}{\nu^2}\right)|x-y|.
$$
By reversing the roles of $x$ and $y$ we get the full inequality  \eqref{ddl}.  

In case (iii) we need the following semi-homogeneity of degree two of $\bar{H}$:
\begin{equation}\label{semi-homo}
\mu \bar{H}(x,\frac{p}{\mu}) \geq \bar{H}(x,\frac{p}{\sqrt{\mu}}) \quad \forall 0<\mu<1, x,p \in \mathbb{R}^n.
\end{equation}
This follows from the representation formula \eqref{eqn:barh1}, because Jensen inequality gives
$$
\mu \log E\left[e^{\int_0^t |\sigma^T(x,Y_s)\frac{p}{\mu}|^2 \, ds} |Y_0=y\right] \geq \log E\left[e^{\mu \int_0^t |\sigma^T(x,Y_s)\frac{p}{\mu}|^2 \, ds} |Y_0=y\right] 
$$
and the conclusion is reached after dividing by $t$ and letting $t \rightarrow \infty$. The other ingredient of the proof is the first inequality in \eqref{est_x_H} that relates the regularity in $x$  of $\bar{H}$ with that of the pseudo-coercive Hamiltonian $|\sigma^T(x,y)p|^2$. With these two inequalities one can repeat the proof of the comparison principle for the pseudo-coercive Hamiltonian by 
 Barles and Perthame, see \cite{BP} for the stationary case and \cite{Bal} for the evolutionary case. 
Let us give a sketch of the main points of the proof. We show that for $\mu<1$, $\mu$ sufficiently near to $1$, it holds
$$
\sup_{\R^n\times [ 0,T]} ( u-\mu v)\leq \sup_{\R^n}(  u-\mu v)(\cdot, 0).
$$ If this is true, then the inequality holds also for  $\mu=1$, proving the Theorem.
By contradiction, we assume that for every $\mu<1$,   there exists $(\overline{x}, \overline{t})$ such that 
\begin{equation}\label{contr}  u(\overline{x}, \overline{t})-\mu v(\overline{x}, \overline{t})>  \sup_{\R^n}( u-\mu v)(\cdot, 0).\end{equation}
Let
\[\Phi(x,z,t,s)=  u(x,t)-\mu v(z,s)-\frac{|x-z|^2}{\epsilon^2}-\frac{|t-s|^2}{\eta^2}-\delta \log(1+|x|^2+|z^2|)+ \alpha \mu s.\]
For $\eps, \eta$ small enough, $\Phi$ has a maximum point, that we denote with  $(x' ,z' ,t' ,s') $. 
By standard arguments,  we get  $ \frac{|x'-z'|^2}{\epsilon^2},\frac{|t'-s'|^2}{\eta^2}\longrightarrow 0$ as $\eps, \eta \to 0$.  

If either $s'=0$ or $t'=0$, it is easy to see that we get a contradiction with \eqref{contr}. 
So we   consider the case $(x',z',t',s')\in \R^n\times\R^n\times \mathopen( 0,T\mathclose)\times \mathopen( 0,T\mathclose)$. 
Let \[p=2\frac{x'-z'}{\epsilon^2}, \quad q_x=  \frac{2x'}{1+|x'|^2+|z'|^2},\quad q_z= \frac{-2y'}{1+|x'|^2+|z'|^2},\quad r=2 \frac{t'-s'}{\eta^2}.\] 
Using the fact that $u$ is a subsolution  we get 
\begin{equation}\label{sub1}
r-  \bar{H}(x',p+\delta q_x)\leq 0.
\end{equation}
Since  $v$ is a supersolution  and  $\bar{H}$ satisfies \eqref{semi-homo}, we get 
\begin{equation}\label{sup}\frac{r}{\mu}-\frac{1}{\mu }\bar{H}\left(z',\frac{p+\delta q_z}{\sqrt{\mu}}\right)\geq   \alpha \end{equation}
So, we multiply \eqref{sup}  by $-\mu$  and sum up to \eqref{sub1} to obtain 
\begin{equation}\label{eqn:disu1}   \bar{H}\left(z',\frac{p+\delta q_z}{\sqrt{\mu}}\right)-\bar{H}(x', p+\delta q_x )  \leq   -\alpha\mu. \end{equation}
Using \eqref{est_x_H} we get  \begin{equation} \label{e1} \bar{H}\left(z', \frac{p+\delta q_z}{\sqrt{\mu}}\right)-\bar{H}(x', p+\delta q_x )  \geq \min_{y \in \mathbb{R}^n}\left(\frac{1}{\mu} |\sigma^T(z',y) (p+\delta q_z)|^2-|\sigma^T(x',y)(p+\delta q_x)|^2\right).\end{equation} 
Let
\[
 A(y)=|\sigma^T(z',y)
(p+\delta  q_z)|\]\[ \Delta(y)= ((\sigma(x',y)-\sigma(z',y))^T(p+\delta q_z),\quad J(y)=\delta \sigma^T(x',y)
(q_x-q_z).
\]
 Note that $\Delta(y)$ goes to zero for $\epsilon,\eta\to 0$ and for all $\delta$ fixed uniformly in $y$, and $J(y)$ goes to zero for $\epsilon,\eta,\delta\to 0$ uniformly in $y$. 
Then we can rewrite the rhs of \eqref{e1} as 
 \begin{equation}\label{eq:1.1.1}\min_{y\in\R^m} \left( \frac{1-\mu}{\mu} A(y)^2 -\Delta(y)^2-J(y)^2-2A(y)\Delta(y)-2J(y) \Delta(y)-2 A(y) J(y)\right ).\end{equation}
  Moreover, for all  $k_1,k_2>0$  and for all $y\in \mathbb{R}^m$ it holds
\[ -2A(y)\Delta(y) \geq -k_1 A(y)^2 -\frac{1}{k_1}\Delta(y)^2,\quad\text{and  } \quad  -2A(y)J(y) \geq -k_2 A(y)^2 -\frac{1}{k_2}J(y)^2.\]
So, recalling \eqref{eqn:disu1}, \eqref{e1} and  \eqref{eq:1.1.1} we get 
\[-\alpha \mu\geq  \min_{y\in\R^m} \left( \left(\frac{1-\mu}{\mu}-k_1-k_2\right) A(y)^2 -\left(1+\frac{1}{k_1}\right)\Delta(y)^2-\left(1+\frac{1}{k_2}\right)J(y)^2-2J(y) \Delta(y)\right).\]
If we choose  $k_1,k_2>0$ such that $k_1+k_2<\frac{1-\mu }{\mu}$ 
then we obtain \[0>-\alpha \mu\geq  \min_{y\in\R^m} \left(  -\left(1+\frac{1}{k_1}\right)\Delta(y)^2-\left(1+\frac{1}{k_2}\right)J(y)^2-2J(y) \Delta(y)\right)\to 0,\]  as  $\eps, \eta, \delta\to 0$,  reaching a contradiction.

\end{proof}

\section{The supercritical case: $\alpha>2$}\label{sec:casosupercritico}

As  in Section \ref{sec:casocritico}, we prove the existence of an effective Hamiltonian giving the limit PDE and first we identify the cell problem that we wish to solve. 
 Plugging the asymptotic expansion  
$$
v^{\eps}(t,x,y)=v^{0}(t,x) + \eps^{\alpha-1} w(t,x,y)
$$
in the equation \eqref{eqn:equazione3}   we get
$$
v^0_t = |\sigma^T D_xv^0|^2 + b \cdot D_yw +\tr(\tau \tau^T D_{yy}^2 w)+O(\eps).
$$
We  consider the \textit{$\delta$-cell problem} for fixed $(\bar{x},\bar{p},\bar{X})$
\begin{equation}\label{eqn:deltacell}
\delta w_{\delta}(y) - |\sigma(\bar{x}, y)^T \bar{p}|^2 - b(y) \cdot D_y w_{\delta}(y)- \tr(\tau(y) \tau(y)^T D_{yy}^2 w_{\delta}(y))=0\, \, \mbox{in} \, \, \mathbb{R}^m,
\end{equation}
where $w_{\delta}$ is the \textit{approximate corrector}. \\
The next result states that $\delta w_{\delta}$ converges to $\bar{H}$ and it is smooth.
\begin{prop}
For any fixed $(\bar{x},\bar{p})$ there exists a constant $\bar{H}(\bar{x},\bar{p})$ such that $\bar{H}(\bar{x},\bar{p})=\lim_{\d\to 0} \d w_\d(y)$ uniformly, where $w_\d \in C^2(\mathbb{R}^m)$ is the unique periodic solution of \eqref{eqn:deltacell}. Moreover
\begin{equation}\label{eqn:effham2}
\bar{H}(\bar{x},\bar{p}) := \int_{\mathbb{T}^m} \! |\sigma(\bar{x}, y) ^T \bar{p}|^2 \, d\mu(y) \, \, \mbox{ uniformly in} \, \, \mathbb{T}^m,
\end{equation}
where $\mu$ is the invariant probability measure on $\mathbb{T}^m$ of the stochastic process \[dY_t=b(Y_t)dt+\sqrt{2}\tau(Y_t)dW_t,\] that is, the periodic solution of  
\begin{equation}\label{eqn:adj}
-\sum_{i,j} \frac{\partial^2}{\partial y_i \partial y_j} ((\tau \tau^T )_{ij}(y) )\mu + \sum_i \frac{\partial}{\partial y_i}( b_i(y) )\mu=0  \qquad  \mbox{in}\, \, \mathbb{R}^m,  
\end{equation}
with  $\int_{\mathbb{T}^n} \! \mu(y) \, dy =1$.
\end{prop}
\begin{proof}
The proof\, essentially follows the arguments presented in \cite{AL, BA2} of ergodic control theory in periodic enviroments.
\end{proof}

\begin{oss}\upshape
Note that in dimension $n=1$ the effective Hamiltonian  assumes the 
 form
\[
H(\bar{x},\bar{p}) = \int_{\mathbb{T}^m}\sigma(\bar{x}, y)^2 d\mu(y)\bar{p}^2= (\bar{\sigma} \bar{p})^2,
\]
where $\bar{\sigma}= \sqrt{\int_{\mathbb{T}^m}  \sigma(\bar{x}, y)^2 d\mu(y)}$. 
\end{oss}

We list some elementary properties of the effective Hamiltonian $\bar{H}$.

\begin{prop}\label{prop:contbarh2} 
 $\bar{H}$ satisfies   properties (a), (b), (c), (d) as  in  Proposition \ref{prop:contbarh1}.
Moreover  
\begin{enumerate}\itemsep2pt
\item[(f)] for all $x, \bar x, p, \bar p\in\R^n$, 
\begin{multline}
\label{est_x_H2}
\min_{y\in\R^m}\left(|\sigma^T(x,y)p|^2 - |\sigma^T(\bar x,y)\bar p|^2\right)\leq \bar H(x,p) - \bar H(\bar x,\bar p)\\ \leq \max_{y\in\R^m}\left(|\sigma^T(x,y)p|^2 - |\sigma^T(\bar x,y)\bar p|^2\right);
\end{multline}
\item[(g)] for every $\lambda\in \R$,  $x,   p \in\R^n$, \begin{equation}\label{2hom}\bar{H}(x,\lambda p)= |\lambda|^2  \bar{H}(x,p).\end{equation}  
\end{enumerate}
\end{prop}
\begin{proof} For the proofs of  (a), (b), (c), (d) we repeat the same arguments as in Proposition \ref{prop:contbarh1}. Properties (f), (g) can be easily checked from the representation formula \eqref{eqn:effham2}.
\end{proof}

We now state the comparison principle among viscosity sub- and supersolutions of the 
 limit PDE 
 \begin{equation}\label{eqn:effPDE2}
 v_t -\int_{\mathbb{T}^m} \! |\sigma(x, y) ^T Dv|^2 \, d\mu(y)=0 \quad \mbox{in} \, \, (0,T) \times \mathbb{R}^n. 
\end{equation}
In this case, differently from the critical  case, we do not need additional assumptions for th e comparison principle to hold.  
\begin{thm}
\label{thm:compbarh2} Let $u\in BUSC([0, T]\times\R^n)$ and $v\in BLSC([0,T]\times \R^n)$ be, respectively, a bounded upper semicontinuous subsolution and   a bounded lower semicontinuous supersolution to \eqref{eqn:effPDE2}
 such that $u(0,x)\leq 
  v(0,x)$ for all $x\in \R^n$. Then $u(x,t)\leq v(x,t)$ for all $x\in \R^n$ and $0\leq t\leq T$. 
\end{thm}
\begin{proof} The homogeneity \eqref{2hom} of the Hamiltonian $\overline{H}$ implies \eqref{semi-homo}, moreover 
\eqref{est_x_H2} holds. 
Then the proof of Theorem \ref{thm:compbarh1}, case (iii), applies here.
\end{proof}

 \section{The subcritical case: $\alpha<2$}\label{5}
 \subsection{The effective Hamiltonian}
 
In this case, the asymptotic expansion we plug in the equation is 
\begin{equation}
\label{eqn:asy}
v^{\eps}(t,x,y)=v^0(t,x) + \eps^{\frac{\alpha}{2}}w(t,x,y).
\end{equation}
Plugging \eqref{eqn:asy} into the equation \eqref{eqn:equazione3} 
we get 
\begin{equation}
\label{eqn:sub1}
v_t^0 =|\sigma^T D_{x}v^0|^2 + 2(\tau \sigma^T D_x v^0 ) \cdot D_y w +  | \tau^T D_{y} w|^2+ O(\eps) .
\end{equation} 
Therefore the cell problem  we want to solve is finding, 
for any fixed $(\bar{x}, \bar{p})$, 
a unique constant $\bar{H}$ such that there is a viscosity solution $w$ of the following equation
\begin{equation}\label{eqn:cell}
\bar{H}(\bar{x}, \bar{p}) -  2(\tau(y) \sigma(\bar{x}, y)^T \bar{p}) \cdot D_y w(y) - |\tau(y)^T D_y w(y)|^2 -|\sigma(\bar{x},y)^T \bar{p}|^2=0.
\end{equation}
Since 
\[2(\tau(y) \sigma^T(\bar{x},y) \bar{p}) \cdot D_y w = 2 (\sigma^T(\bar{x},y) \bar{p})\cdot (\tau^T(y) D_y w) \],
we can restate the cell problem as 
\begin{equation}\label{eqn:cellnew}
\bar{H}(\bar{x}, \bar{p}) - |\tau^T(y) D_y w(y)+\sigma^T(\bar{x},y)  \bar{p}|^2=0 .
\end{equation}

The following proposition deals with the existence and uniqueness of  $\bar{H}$.
\begin{prop}\label{prop51}
For any fixed  $(\bar{x},\bar{p})$, there exists a unique constant $\bar{H}(\bar{x},\bar{p})$ 
such that the  cell problem \eqref{eqn:cell} admits a   periodic viscosity  solution $w$. Moreover $w$ is  Lipschitz continuous and 
 there exists  $C>0$ independent of $\bar{x}, \bar{p}$ such that 
$$
\max_y|Dw(y;\bar{x},\bar{p})|\leq C(1+|\bar p|).
$$
\end{prop}
\begin{proof} As for the other cases  
we introduce the following approximant problem, with $\delta >0$, 
\begin{equation}
\label{eqn:cellappro}
\delta w_{\delta}(y) - |\tau^T(y)  D_y w_\delta (y)+\sigma^T(\bar{x},y) \bar{p}|^2=0 
\,\mbox{ in } \mathbb{R}^m.
\end{equation} 
Let $w_\d$ the unique periodic viscosity solution to \eqref{eqn:cellappro}. By standard comparison principle we get that 
$$
|\d w_\d|\leq \max_{y\in\R^m}|\sigma^T(\bar{x},y)\bar{p}|^2\leq C(1+|\bar{p}|^2)\quad \forall y \in \mathbb{R}^m.
$$ 
Moreover, using the coercivity of the Hamiltonian (see \cite[Prop II.4.1]{BCD}), we get that $w_\d$ is Lipschitz continuous and there exists a constant $C$ independent of $\d$ and $\bar{p}$ such that 
$$
\max_{y\in\mathbb{R}^m}|Dw_\d |\leq C(1+|\bar{p}|).
$$  
So, we conclude as in the proof of Proposition \ref{thm:trucell}.
 \end{proof} 
We give some representation formulas for the effective Hamiltonian $\bar{H}$. 
\begin{prop} 	
(i) $\bar{H}$ satisfies 
\begin{equation}
\label{eqn:barhsub}
\bar{H}(\bar{x}, \bar{p})= \lim_{\delta \rightarrow 0} 
\sup_{\beta(\cdot)
}  \delta  \int_0^{+\infty} \! \left(|\sigma(\bar{x},y
(t))^T\bar{p}|^2 - |\beta(t)|^2\right)e^{-\delta t} \, dt,
\end{equation}
where $\beta(\cdot)$ varies over measurable functions taking values in $\R^r$, $y
(\cdot)$ is the trajectory of the 
control system
\begin{equation}
\left\{
\begin{array}{ll}
\dot{y}(t)= 2 \tau(y(t)) \sigma^T(\bar{x}, y(t)) \bar{p}-2\tau(y(t)) \beta , \quad t>0 ,\nonumber\\
 y(0)=y 
\end{array}
\right.\,
\end{equation}
and the limit is uniform with respect to the initial position $y$ of the system.

(ii) If, in addition, $\tau(y) \sigma^T({x}, y)=0$ for all $x, y$, then
 \begin{equation}
\label{eqn:barhunc}
\bar{H}(\bar{x}, \bar{p})= \max_{y\in\R^m}|\sigma^T(\bar{x},y)\bar{p}|^2 .\end{equation}

(iii) If  $n=m=r=1$,  and $\sigma\geq 0$
\begin{equation}\label{eqn:lpv}  
\bar{H}(\bar{x}, \bar{p})= 
\left(\int_0^1 \frac{\sigma(\bar{x},y)}{\tau(y)} \,dy\right)^2\left(\int_0^1 \frac{1}{\tau(y)} \,dy\right)^{-2} \bar{p}^2.
\end{equation}

\end{prop}
\begin{proof}

 The formula \eqref{eqn:barhsub} can be 
 proved by writing 
 \eqref{eqn:cellappro} as a Bellman equation
\begin{equation}\label{previous}
\delta w_\delta(y)  + \inf_{\beta \in \mathbb{R}^r} \left\{ \left(2\tau(y) \beta -2\tau(y)\sigma(\bar{x},y)^T \bar{p}\right)\cdot D_y w_\delta+|\beta|^2\right\} 
- |\sigma(\bar{x},y)^T \bar{p}|^2=0.
\end{equation}
Then $w_\delta$ is the value function of the infinite horizon discounted deterministic control problem appearing in \eqref{eqn:barhsub} (see, e.g., \cite{BCD, Barl}).

If $\tau(y) \sigma^T({x}, y)=0$ for all $x, y$, then \eqref{eqn:cellnew} reads 
\[
- |\tau^T(y) D_y w(y)|^2= |\sigma^T(\bar{x},y)  \bar{p}|^2 -\bar{H}(\bar{x},\bar{p}).
\]
So, this gives immediately the  inequality  $\geq$ in \eqref{eqn:barhunc}. 
The other inequality is obtained by standard comparison principle arguments applied to 
the approximating problem \eqref{eqn:cellappro}. 

Finally, in the case $n=m=r=1$,
if $\bar{p}\geq 0$ we write explicitly the corrector as 
$$
w(y)=\int_{0}^y \! \frac{\bar{H}^{\frac{1}{2}} -\sigma(\bar{x},s) \bar{p} }{\tau(s)} \, ds . \,
$$
Note that $w\in C^1$ is periodic and does the job. A similar construction works for $\bar{p}<0$.

\end{proof} 

For the comparison principle it is useful to define 
\[
H_0(\bar{x}, \bar{p})=\sqrt{\bar{H} (\bar{x}, \bar{p})}
\]
and 
observe that the cell problem \eqref{eqn:cell} is equivalent to the following equation
\begin{equation}\label{eqn:cellsq}
H_0(\bar{x}, \bar{p}) - |\tau^T(y)  D_y w(y)+\sigma^T(\bar{x},y)  \bar{p}|=0 .
\end{equation}
Here are some properties of 
 $\bar{H}$ and $H_0$.
\begin{prop}\label{prop:contbarh3} 
 $\bar{H}$ satisfies   properties (a), (b), (c), (d) as  in  Proposition \ref{prop:contbarh1}.
Moreover 
$\bar{H}(x,p)=(H_0(x,p))^2$ 
with $H_0$ positively $1$ homogeneous, i.e.
 \begin{equation}\label{ho} H_0(x,\lambda p)=|\lambda| H_0(x,p)\qquad \forall\lambda \in\R ,
 \end{equation} 
 there exists $C>0$ such that  $|H_0(x,p)|\leq C|p|$, and
\begin{equation}\label{li} |H_0(x,p)-H_0(z,p)|\leq C(1+|p|)|x-z|\qquad \forall x,z\in\R^n, \ p\in\R^n.\end{equation} 
 \end{prop}
\begin{proof} 
For the proofs of  (a), (b), (c) we repeat the same arguments as in Proposition \ref{prop:contbarh1}. The properties of $H_0$ defined in  \eqref{eqn:cellsq} follow from standard theory, using comparison type argument in the approximating problem 
\[
\delta v_{\delta}(y) - |\tau^T(y)  D_y v_\delta (y)+\sigma^T(x,y) p|=0 
\,\mbox{ in } \mathbb{R}^m.\] 
 \end{proof}

\subsection{Comparison principle} 
We consider the limit PDE 
 \begin{equation}\label{eqn:effPDE3}
 v_t - \bar{H}(x,Dv)=0 \quad \mbox{in} \, \, (0,T) \times \mathbb{R}^n. 
\end{equation}
We now state the comparison principle for the effective Hamiltonian $\bar{H}$. 
\begin{thm}
\label{thm:compbarh3} Let $u\in BUSC([0, T]\times\R^n)$ and $v\in BLSC([0,T]\times \R^n)$ be, respectively, a bounded upper semicontinuous subsolution and   a bounded lower semicontinuous supersolution to \eqref{eqn:effPDE}
 such that $u(0,x)\leq h(x)\leq v(0,x)$ for all $x\in \R^n$. Then $u(x,t)\leq v(x,t)$ for all $x\in \R^n$ and $0\leq t\leq T$. 
\end{thm}
\begin{proof} 
Recall that  $\bar{H}=H_0^2$ and $H_0$ is continuous and satisfies \eqref{ho} and \eqref{li}. So, we can apply Theorem 2.4 in \cite{cdl}. 
\end{proof}

\section{The convergence result}\label{6}
In this Section we state  the main result of the paper, namely, the convergence theorem for the singular perturbation problem. We will make use of the 
 relaxed semi-limits which we define as follows.
For the functions $v_\eps$ introduced in Section \ref{log-tran} 
 the relaxed upper  semi-limit $\bar{v}=\limsup^*_{\eps \rightarrow 0}\sup_{y} v^{\eps}$ is 
$$
\bar{v}(t,x) := \limsup_{\eps \rightarrow 0, (t^{'}, x^{'}) \rightarrow (t,x)} \sup_{y} v^{\eps}(t^{'}, x^{'}, y) , \quad 
x\in\R^n , \,t\geq 0 .
$$
We define analogously the lower semi-limit $\underline{v}=\liminf_{*\eps \rightarrow 0}\inf_{y} v^{\eps}$ by replacing $\limsup$ with $\liminf$ and $\sup$ with $\inf$. 
Since $h$ is bounded the family $v^\eps$ is equibounded and we have $\bar{v}\in BUSC([0, T]\times\R^n)$ and $\underline{v}\in BLSC([0, T]\times\R^n)$.

The standing hypotheses of sections \ref{asssistema} and \ref{log-tran} are assumed in this section.
\subsection{The convergence result: critical and supercritical case, $\alpha\geq 2$.}
Recall that by Proposition \ref{prop:propequazione} i) $v^{\eps}$ defined by \eqref{v-eps} is the solution of   
\begin{equation*}\label{eqn:eq1} \begin{cases} 
\partial_t v^\eps  - H^\eps\left(x,y, D_x v^\eps, \frac{D_y v^\eps}{\eps^{\alpha-1}}, D_{xx} v^\eps, \frac{D_{yy}^2 v^\eps}{\eps^{\alpha-1}}, \frac{D_{xy} v^\eps}{\eps^{\frac{\alpha-1}{2}}}\right)=0  &   (0, T)\times\R^n\times\R^m\\ v^\eps(0,x,y)=h(x) &   \R^n\times\R^m.\end{cases}
\end{equation*}
with
\begin{eqnarray*} 
H^{\eps}(x,y, p, q, X,Y, Z):&=& |\sigma^T p|^2+ b\cdot q+ \tr(\tau \tau^T Y) + \eps\left(\tr(\sigma\sigma^T X) + 
\phi \cdot p\right) \\ &+& 2\eps^{\frac{\alpha}{2}-1} (\tau\sigma^ T p ) \cdot q +2\eps^{\frac{1}{2}}\tr(\sigma\tau^T Z) + \eps^{\alpha-2}|\tau^T q|^2 .
\end{eqnarray*}
\begin{thm}\label{thm:teoconv1}
Assume $\alpha\geq 2$. Then

i) The upper limit $\bar{v}$ (resp., the lower limit $\underline{v}$) of $v^{\eps}$  is a subsolution (resp., supersolution) of the effective equation
\begin{equation}\label{eqn:effprob1}
v_{t} - \overline{H}(x, Dv)=0 \, \, \mbox{in} \, \, (0,T) \times \mathbb{R}^n \quad v(0,x)=h(x) \, \, \mbox{on} \, \, \mathbb{R}^n
\end{equation}
 where $\bar{H}$ is given by \eqref{eqn:effham2} for $\alpha >2$, and it is defined by Proposition \ref{thm:trucell} for $\alpha=2$ (with the formulas \eqref{eqn:khaise}, \eqref{khaiset}, \eqref{meas}, and \eqref{eqn:barh1});
 
 ii) if $\alpha >2$ then $v^\eps$ converges uniformly on the compact subsets of $[0,T) \times \mathbb{R}^n \times \mathbb{R}^m$ to the unique viscosity solution of \eqref{eqn:effprob1}.
 
 iii) if $\alpha =2$ and 
\begin{equation}\label{additional1}
\begin{cases}\text{either  $\sigma=\sigma(y)$ is independent of $x$ and $h\in BUC(\R^n)$,}\\
\text{or, for some $\nu>0$,  $|\sigma^T(x,y)p|^2>\nu |p|^2 \quad\forall \, x, p\in\R^n, \,y\in\R^m$},\\
\text{or, $\tau(y)\sigma^T(x,y)=0$ for all $x,y$},
\end{cases}
\end{equation}
 then $v^\eps$ converges uniformly as in ii).
\end{thm}
 
\begin{proof}  
{\it i)} The inequalities $\underline{v}(0,x)\leq h(x)\leq\bar{v}(0,x)$ follow from the definitions. 
The problem of taking the limit in the PDE is a regular perturbation of a singular perturbation problem, in the terminology of  \cite{ABM}.  The result can be 
proved  by the methods developed in \cite{ABM} for such problems, with minor modifications.
%

{\it ii)} By the definition of the semilimits 
$
\underline{v} \leq \bar{v} 
$
 \,\,in $[0,T) \times \mathbb{R}^n$. 
The comparison principle Proposition \ref{thm:compbarh2} for the effective equation \eqref{eqn:effprob1} gives the inequality $\leq$ 
and therefore 
$
\bar{v}=\underline{v}=v 
$ 
 in $[0,T] \times \mathbb{R}^n$. Thanks to the properties of semilimits, we finally get that $v^\eps$ converges locally uniformly to the 
unique bounded solution of \eqref{eqn:effprob1}. 

{\it iii)} The proof is the same as for {\it ii)}, but now we need the additional assumption \eqref{additional1} for the comparison principle Theorem \ref{thm:compbarh1}.
\end{proof}

\subsection{The convergence result: subcritical case, $\alpha<2$.}
Recall that by Proposition \ref{prop:propequazione} ii) $v^{\eps}$ defined by \eqref{v-eps} is the solution of   
\begin{equation*}
 \begin{cases} 
v^{\eps}_t  =H_\eps\left(x,y, D_x v^\eps, \frac{D_y v^\eps}{\eps^{\frac{\alpha}{2}}}, D_{xx} v^\eps, \frac{D_{yy}^2 v^\eps}{\eps^{\frac{\alpha}{2}}}, \frac{D_{xy} v^\eps}{\eps^{\frac{\alpha}{4}}}\right)&   (0, T)\times\R^n\times\R^m\\ v^\eps(0,x,y)=h(x) &   \R^n\times\R^m.
\end{cases}
\end{equation*}
with
\begin{eqnarray}\nonumber  
H_{\eps}(x,y, p, q, X,Y, Z):&=& |\sigma^T p|^2+ 2
(\tau \sigma^Tp)\cdot q+  |\tau^T q|^2 +\eps\left(\tr(\sigma\sigma^T X) + 
\phi \cdot p\right) \\ \nonumber &+&   2\eps^{1-\frac{\alpha}{4}}\tr(\sigma\tau^T Z)  + \eps^{1-\frac{\alpha}{2}}b\cdot q +  \eps^{1-\frac{\alpha}{2}}\tr(\tau \tau^T Y)\nonumber. 
\end{eqnarray}
\begin{thm}\label{thm:teoconv2}
Assume $\alpha<2$. Then

i) the upper limit $\bar{v}$ (resp., the lower limit $\underline{v}$) of $v^{\eps}$  is a subsolution (resp., supersolution) of the effective equation
\eqref{eqn:effprob1}
 where $\bar{H}$ is defined by Proposition \ref{prop51}  
(with the formula \eqref{eqn:barhsub});

ii) $v^\eps$ converges uniformly on the compact subsets of $[0,T) \times \mathbb{R}^n \times \mathbb{R}^m$ to the unique viscosity solution of \eqref{eqn:effprob1}.
\end{thm}
 \begin{proof}The proof is the same as that of  Theorem \ref{thm:teoconv1}, by using the comparison principle Proposition \ref{thm:compbarh3}.
\end{proof}

 
\begin{oss}\rm{
In the case $\alpha \leq 2$  we can give a convergence result analogous to Theorem \ref{thm:teoconv1} and Theorem \ref{thm:teoconv2} 
for a terminal cost $h=h(x,y)$ depending also on the fast variable $y$, so that
the payoffs is
\begin{equation}\label{vepsy}
v_{\eps}(t,x,y) := \eps \log E\left[e^{\frac{h(X_{t}, Y_t)}{\epsilon}} | (X.,Y.)\, \, \mbox{satisfy \eqref{sistema0}}\right],
\end{equation}
In this case 
we must find a suitable \textit{effective initial value} $\bar{h}$ depending only on the variable $x$; moreover  the convergence cannot be up to time $t=0$
but only on the compact subsets of $(0,T) \times \mathbb{R}^n \times \mathbb{R}^m$
 to the unique viscosity solution of 
$$
v_{t} - \overline{H}(x, Dv) =0\, \, \mbox{in} \, \, (0,T) \times \mathbb{R}^n \quad v(0,x)=\bar{h}(x) \, \, \mbox{on} \, \, \mathbb{R}^n.
$$
The proof follows the methods of \cite{BA}, where an asymptotic problem for finding $\bar{h}$ is given and the 
relaxed semi-limits are modified at $t=0$ 
 to 
 deal with the expected initial layer. For further details and proofs we refer to \cite{TD}.}
\end{oss}


\section{The large deviation principle}\label{7}

In this section we derive a large deviation principle for the process $X^\eps_t$ defined in \eqref{eqn:systemscaled}. Throughout the section we suppose that $\sigma$ is uniformly non degenerate, that is,
for some $\nu>0$ and for all $x,p \in \mathbb{R}^n$
\begin{equation}\label{coercive2}
|\sigma^T(x,y)p|^2>\nu|p|^2.
\end{equation}
By \eqref{eqn:bounds}, under \eqref{coercive2}, the effective Hamiltonian is coercive. 
Let $\bar{L}$ be the \textit{effective Lagrangian}, i.e. for $x \in \mathbb{R}^n$
\begin{equation}\label{lagr}
\bar{L}(x,q)=\max_{p\in \mathbb{R}^n} \{p\cdot q-\bar{H}(x,p)\}.
\end{equation}
Note that $\bar{L}(x,\cdot)$  is a convex nonnegative function such that $\bar{L}(x,0)=0$ for all $x \in \mathbb{R}^n$, since $\bar{H}(x,\cdot)$ is  convex nonnegative and $\bar{H}(x,0)=0$ for all $x \in \mathbb{R}$.\\
For each $x_0 \in \mathbb{R}^n$ and $t>0$, define
\begin{equation}\label{I}
I(x;x_0,t):= \inf\left[\int_0^t \! \bar{L}\left(\xi(s),\dot{\xi}(s)\right)\, ds \ \Big|\ \xi\in AC(0,t), \ \xi(0)=x_0, \xi(t)=x\right] .
\end{equation}

\begin{oss}\label{oss:rate}\rm{
(a) The function $I$ defined in \eqref{I}  is continuous in the variable $x$ (see, e.g., \cite{MF}) 
and is a nonnegative function such that $I(x_0;x_0,t)=0$. 

(b) $I$ satisfies the following growth condition for some $C>0$ and all $x, x_0 \in \mathbb{R}^n$
\begin{equation}\label{eqn:growth}
\frac{1}{4C} \frac{|x-x_0|^2}{t} \leq I(x;x_0,t)\leq \frac{1}{4\nu} \frac{|x-x_0|^2}{t},
\end{equation}
where $\nu$ is defined in \eqref{coercive2}. In fact, thanks to the property \eqref{eqn:bounds} stated in Proposition \ref{prop:contbarh1}, we get that
$$
\frac{1}{4C}|p|^2\leq \bar{L}( x, p)\\ \leq \frac{1}{4 \nu}|p|^2.
$$
Then we have
$$
\frac{1}{4C}\inf_{\xi(0)=x_0, \xi(t)=x } \int_0^t \!  |\dot{\xi}(s)|^2\leq I(x;x_0,t)\leq\frac{1}{4\nu}\inf_{\xi(0)=x_0, \xi(t)=x } \int_0^t \!  |\dot{\xi}(s)|^2,
$$
from which we get \eqref{eqn:growth}.

(c) If $\sigma$ does not depend on $x$, i.e.  $\bar{H}=\bar{H}(p)$,  the rate function in \eqref{I} is 
$$
I(x;x_0,t)=t\bar{L}\left(\frac{x-x_0}{t}\right).
$$

(d) If $\sigma$ does not depend of $x$ and $n=1$,  $I$ is a  monotone nondecreasing function of $x$ when $x > x_0$. Analogously, $I$ is a monotone nonincreasing function of $x$ when $x < x_0$. 
}
\end{oss}

\begin{thm}\label{thm:ldp}
Let $(X^\eps, Y^\eps)$ be the process defined in \eqref{eqn:systemscaled} with initial position $X^\eps_0=x_0$ and $Y^\eps_0=y_0$. 
Then for every  $t>0$, a large deviation principle holds for $\{X^\eps_t:\eps >0\}$ with speed $\frac{1}{\eps}$ and good rate function $I(x;x_0, t)$.  In particular, for any open set $B \subseteq \mathbb{R}^n$
\begin{equation}\label{ldp}
\lim_{\epsilon \rightarrow 0} \epsilon \log P(X^\eps_t \in B)=-\inf_{x \in B} I(x;x_0,t).
\end{equation}
\end{thm}

\begin{oss} \upshape \label{mono} 
 Thanks to Remark \ref{oss:rate}, if $\sigma$ does not depend on $x$ and $n=1$, we have $\inf_{y>x}I(y;x_0,t)=I(x;x_0,t)$ for $x\geq x_0$ and \eqref{ldp} can be written in the following way
$$
\lim_{\eps \rightarrow 0} \eps \log P(X^\eps_t > x) = -I(x; x_0, t) \quad \mbox{when} \, \, x > x_0
$$ 
and analogously when  $x < x_0$
$$
\lim_{\eps \rightarrow 0} \eps \log P(X^\eps_t < x) = -I(x; x_0, t).
$$ 
\end{oss}

\begin{oss}\rm{
We note that the rate function $I$ defined in \eqref{I} does not depend on the drift $\phi$ of the log-price $X^{\eps}_t$ and it depends only on the volatility $\sigma$ and on the fast process $Y^{\eps}_t$. In fact, this holds for the effective Hamiltonian $\bar{H}$ by 
the representation formulas 
 \eqref{eqn:khaise} for $\alpha=2$,  \eqref{eqn:effham2} for  $\alpha>2$ and \eqref{eqn:barhsub} for $\alpha<2$, 
 and hence it holds for the Legendre transform $\bar{L}$. }
\end{oss}

\begin{proof}
 We divide the proof  in two steps, the first is the proof of the large deviation principle, while the second is the proof of the representation formula \eqref{I} for the good rate function.  

\begin{step}1 (Large deviation principle) \upshape The proof of this step is  similar to that of Theorem $2.1$ of \cite{FFK} with some minor changes.
The idea is to apply Bryc's inverse Varadhan lemma (see Appendix A, Lemma \ref{lem:bryc}) with  $\mu_\eps$   given by the laws of $\{X^\eps_t\}$ and $\Lambda_h^\eps$ given by $v_\eps$. Recall that, for $h \in BC(\mathbb{R}^n)$, $v_\eps$ is defined as 
$$
v_{\eps}(t,x,y) := \eps \log E\left[e^{\frac{h(X^\eps_{t})}{\epsilon}} | (X^\eps_.,Y^\eps_.)\, \, \mbox{satisfy \eqref{eqn:systemscaled}}\right] .
$$
We proved in Theorems \ref{thm:teoconv2}, \ref{thm:teoconv2} that $v_\eps$ converge uniformly  to a function $v^h$.  \\
To apply Lemma \ref{lem:bryc},  we have to prove the exponential tightness of $\{X^\eps_t\}$.
 Define the following function
\begin{equation}
f_\eps(x,y)=
\left\{
 \begin{array}{ll}
  f(x) + \eps^{\alpha-1} \zeta(y)\quad &\mbox{if} \, \, \alpha \geq 2,  \\
  f(x) + \eps^{\frac{\alpha}{2}} \zeta(y)\quad &\mbox{if} \, \, \alpha <2,
                 \end{array}
\right.\,
\end{equation}
where
$$
f(x)=\log(1+|x|^2)
$$
and $\zeta(y)$ is a positive differentiable function with bounded first and second derivatives. Since $f(x)$ is an increasing function of $|x|$ and since $\zeta(y) \geq0$, we have that for any $c>0$ there exists a compact set $K_c \subset \mathbb{R}^n$ such that 
\begin{equation}\label{eqn:fc}
f_\eps(x,y) >c\, \, \mbox{ when }x \not \in K_c.
\end{equation}
 We observe that $||\partial_{x_j} f||_{\infty} +||\partial_{x_{j}x_{i}}^2 f||_{\infty} < \infty$ for all $i=1\cdots n, j=1 \cdots n$, and by our choice of $\zeta$ we therefore have that
\begin{equation}\label{suph}
\sup_{x\in \mathbb{R}^n, y \in \mathbb{R}^m} H_\eps (x,y, D_x f_\eps, D_y f_\eps, D^2_{xx} f_\eps, D^2_{yy} f_\eps, D^2_{xy} f_\eps)=C < \infty,
\end{equation}
where $H_\eps$ is defined as follows
\begin{eqnarray*}
H_\eps(x,y,p, q , X, Y, Z) &= &|\sigma^T p|^2+ \eps \tr(\sigma\sigma^T X) + \eps\phi \cdot p + 2\eps^{-\frac{\alpha}{2}}\tr (\tau \sigma^ T p ) \cdot q \\ &+& 2\eps^{1-\frac{\alpha}{2}}\tr(\sigma\tau^T Z) + \eps^{1-\alpha}b\cdot q + \eps^{-\alpha}|\tau^T q|^2 + \eps^{1-\alpha}\tr(\tau \tau^T Y).  
\end{eqnarray*}
We will write $H_\epsilon f_\epsilon(x,y)$  to denote $H_\eps (x,y, D_x f_\eps, D_y f_\eps, D^2_{xx} f_\eps, D^2_{yy} f_\eps, D^2_{xy} f_\eps)$.
The $P$ and $E$ in the following proof denote probability and expectation conditioned on $(X,Y)$ starting at $(x,y)$. Define the  process
\begin{equation}\label{eqn:m}
M^\eps_t=\exp\left\{\frac{f_\eps(X^\eps_t,Y^\eps_t)}{\eps} - \frac{f_\eps(x,y)}{\eps} - \frac{1}{\epsilon}\int_0^t H_\eps f_\eps(X^\eps_s, Y^\eps_s)\, ds \right\}.
\end{equation} 
Then $M_{\eps,t}$ is a supermartingale and hence we can apply the optional sampling theorem (see Appendix \ref{sec:prob}, Theorem \ref{thm:opt}), that is
\begin{equation}\label{eqn:sampl}
1 \geq E\left[M^\eps_t\right].
\end{equation}
Then 
\begin{eqnarray}\label{final}
1\geq E\left[M^\eps_t\, |\, X^\eps_t \notin K_c \right] &\geq& E\left[e^{\frac{(c-f_\eps(x,y)-tC)}{\eps}}\, |\,X^{\eps}_{t} \notin K_c \right]\\ \nonumber &=&P(X^\eps_t \not\in K_c) e^{\frac{(c-f_\eps(x,y)-tC)}{\eps}},
\end{eqnarray}
where we have used \eqref{eqn:fc} and \eqref{suph} to estimate the first and third term in $M^\eps_t$.
Then  we get
$$
\eps \log P(X^\eps_t \not\in K_c) \leq t C + f_\eps(x,y) -c \leq \, \, \mbox{const}\,\,-c
$$
and this finally gives us the exponential tightness of $X^\eps_t$.\\
So, by Bryc's inverse Varadhan lemma (see Appendix \ref{sec:prob}, Lemma \ref{lem:bryc}),  the measures associated to the process $X^\eps_t$ satisfy the LDP with the good rate function
\begin{equation}\label{eqn:rate1}
I(x;x_0,t)=\sup_{h \in BC(\mathbb{R}^n)} \{h(x) -v^h(t,x_0)\}
\end{equation}
and
$$
v^h(t,x_0)=\sup_{x \in \mathbb{R}^n}\{h(x)-I(x;x_0,t)\}.
$$
\end{step}
\begin{step}2 (Representation formula for   the good rate function) \upshape 
The solution $v^h$ to the effective equation
\begin{equation}\label{efft}
\left\{
 \begin{array}{ll}
  v_t- \bar{H}(x,Dv)=0 \quad &\mbox{in} \, \, (0,T) \times \mathbb{R}^n\\
v(0,x)=h(x) \quad &\mbox{in} \, \, \mathbb{R}^n
                 \end{array}
\right.\,
\end{equation}
can be represented through the following formula
\begin{multline}\label{eqn:formula}
v^h(t,x)=\\ \shoveleft{\sup \left\{h(y)-\int_0^t \, \bar{L}\left(\xi(s), \dot{\xi}(s)\right)\, ds\ |\ y\in\R^n, \xi\in AC(0,t), \xi(0)=x, \xi(t)=y\right\},}
\end{multline}
where  $\bar{L}$ is the effective Lagrangian defined in \eqref{lagr}. We refer to \cite{MF} where it is shown that $v^h$ is continuous and is the solution of \eqref{efft}.
We define
\begin{equation}\label{eqn:r}
r(x;x_0,t)=\inf_{\xi(0)=x_0, \xi(t)=x} \int_0^t \! \bar{L}\left(\xi(s),\dot{\xi}(s)\right)\, ds
\end{equation}
Thanks to \eqref{eqn:rate1} and \eqref{eqn:formula}, we can write
\begin{multline}\label{i}
I(x;x_0,t)=\\\shoveleft{r(x;x_0,t) + \sup_{h \in BC(\mathbb{R})} \inf \left \{ h(x) -h(y) + \int_0^t \, \bar{L}\left(\xi(s), \dot{\xi}(s)\right)\, ds - r(x;x_0,t)\right\},}
\end{multline}
where the infimum is over $y \in \mathbb{R}^n$ and absolutely continuous functions $\xi$ such that $\xi(0)=x_0, \xi(t)=y$. Then
$$
I(x;x_0,t) =r(x;x_0,t) + J(x;x_0,t),
$$
where  $J(x;x_0,t):=\sup_{h \in BC(\mathbb{R})} J_h(x;x_0,t)$ and 
$$
J_h (x;x_0,t)=  \inf\left \{ h(x) -h(y) + \int_0^t \, \bar{L}\left(\xi(s), \dot{\xi}(s)\right)\, ds - r(x;x_0,t)\right\}.
$$ 
Taking $y=x$, we obtain $J_h(x;x_0,t) \leq 0$ and therefore $J(x;x_0,t) \leq 0$.
Now we define a function $h_{*} \in BC(\mathbb{R})$ as follows:
$$
h_{*}(y)= r(y;x_0,t) \wedge r(x;x_0,t).
$$
We claim that $h_{*}$ is continuous. 
Then
$ 
J_{h_{*}}(x;x_0,t)=0
$ 
and therefore
$ 
J (x;x_0,t)= 0.
$ 
In conclusion
$$
I(x;x_0,t)=\inf_{\xi(0)=x_0, \xi(t)=x} \int_0^t \! \bar{L}\left(\xi(s),\dot{\xi}(s)\right)\, ds.
$$
Finally, the claim follows from the continuity of the function $r(y;x_0,t)$ in the variable $y$, that can be found, e.g., in  \cite{MF},  Section $4$, Proposition $3.1$ and Corollary $3.4$.
\end{step} \end{proof}

\section{Out-of-the-money option pricing and asymptotic implied volatility}\label{8}

\subsection{Option price}\label{price}
In this section, we give some applications of Theorem \ref{thm:ldp} in dimension $1$ to out-of-the-money option pricing.  In particular, in Corollary \ref{cor:vol1}, we state an asymptotic estimate for the behaviour of the price of out-of-the-money European call option with strike price $K$ and short maturity time $T=\eps t$. 

Let $S^\eps_t$ be the asset price, evolving according to the following stochastic differential system 
\vspace{0.5cm}
\begin{equation}
\left\{
\begin{array}{ll}
d S^\eps_t =  \eps\xi(S^\eps_t, Y^\eps_t)S^\eps_t dt + \sqrt{2\eps}\zeta(S^\eps_t, Y^\eps_t) S^\eps_t dW_t \quad &S^\eps_0=S_0 \in \mathbb{R}_{+}\\
dY^\eps_t= \eps^{1-\alpha} b(Y^\eps_t) dt + \sqrt{ 2 \eps^{1-\alpha}}\tau(Y^\eps_t) dW_t \quad &Y^\eps_0=y_0 \in \mathbb{R}^m,
\end{array}
\right.\,
\end{equation}
where $\alpha >1$, $\tau, b$ are as in \eqref{eqn:systemscaled} and $\xi:\R_+\times\R^m\to \R$, $\zeta:\R_+\times\R^m\to \Mi^{1,r}$ are Lipschitz continuous bounded functions, periodic in $y$.  Observe that $S^\eps_t >0$ almost surely if $S_0>0$.
We define   $X_t^\eps=\log S^\eps_t$. Then $(X^\eps_t, Y^\eps_t) $ satisfies \eqref{eqn:systemscaled}
with \[\phi(x,y)=\xi(e^{x},y)-\zeta (e^x,y) \zeta^T(e^x,y)\qquad \sigma (x,y)=\zeta(e^x,y).\] 
We consider out-of-the-money call option by taking
\begin{equation}\label{eqn:assdato}
S_0 < K \quad \mbox{or} \quad x_0 < \log K.
\end{equation}
Following the argument used in \cite{FFK}, we can derive an option price estimates stated in Corollary \ref{cor:vol1}.
Similarly, by considering out-of-the-money put options, one can obtain the same formula for  $S_0 > K$.
\begin{cor}\label{cor:vol1}
Suppose that  $S_0 < K$. Then, for fixed $t>0$
\begin{equation}\label{price1}
\lim_{\eps \rightarrow 0^{+}} \eps \log E\left[ \left( S^\eps_t -K\right)^{+}\right]=-\inf_{y >\log K} I\left (y;x_0,t\right).
\end{equation}
\end{cor}
\subsection{Implied volatility}
We give an asymptotic estimate of the Black-Scholes implied volatility for out-of-the-money European call option, with strike price $K$, which we denote by $\sigma_{\eps}(t,\log K,x_0) $. \\
We recall that given an observed European call option price  for a contract with strike price $K$ and expiration date $T$, the \textit{implied volatility} $\sigma$ is defined to be the value of the volatility parameter that must go into the Black-Scholes formula to match the observed price.\\
By  arguments similar to those of the ones used in \cite{FFK}, we get the following asymptotic formula.
\begin{cor}\label{cor:vol}
\begin{equation}\label{eqn:vol}
\lim_{\eps \rightarrow 0^{+}} \sigma^2_\eps (t,\log K,x_0)=\frac{(\log K - x_0)^2}{2 \inf_{y > \log K}I(y;x_0,t)t}.
\end{equation}
\end{cor}
Note that the infimum in the right-hand side of \eqref{eqn:vol}, is always positive by assumption \eqref{eqn:assdato} and by \eqref{eqn:growth}.

\begin{oss}\upshape 
When $\zeta(s,y)=\zeta(s)$, then thanks to Remark \ref{mono},   \eqref{price1} simplifies to
$$
\lim_{\eps \rightarrow 0^{+}} \eps \log E\left[ \left( S^\eps_t -K\right)^{+}\right]=- I\left (\log K;x_0,t\right)
$$ 
and \eqref{eqn:vol} reads
$$
\lim_{\eps \rightarrow 0^{+}} \sigma^2_\eps (t,\log K,x_0)=\frac{(\log K - x_0)^2}{2 I(\log K;x_0,t)t}.
$$
\end{oss}

\begin{proof}
By the definition of implied volatility
\begin{eqnarray}\label{eqn:seisei}
 E\left[(S^\eps_t -K)^{+}\right]&=& e^{r\eps t} S_0 \Phi\left(\frac{x_0-\log K+r \eps t + \sigma_\eps^2 \frac{\eps t}{2}}{\sigma_\eps \sqrt{\eps t}} \right) \\ &-& K \Phi\left(\frac{x_0-\log K+r \eps t -\sigma_\eps^2 \frac{\eps t}{2}}{\sigma_\eps \sqrt{\eps t}}\right)\nonumber,
\end{eqnarray}
where $\Phi$ is the Gaussian cumulative distribution function.
Then the proof follows as in \cite{FFK}, using \eqref{eqn:seisei} and Corollary \ref{cor:vol1}.
\end{proof}

\appendix
\section{}\label{sec:prob}

We recall some standard notions from large deviation theory that we need in section \ref{7}.
Throughout the section, $\mu_\epsilon$ will denote a family of probability measures defined on  $\mathbb{R}^n$ with its  Borel $\sigma$-field $\mathcal{B}$. For the definitions and theorems in a more general setting and for further details we refer to \cite{DZ}.\\
Given a family of probability measures $\{\mu_\epsilon\}$, a large deviation principle characterizes the limiting behavior, as $\epsilon \rightarrow 0$, of $\{\mu_\epsilon\} $ in terms of a rate function through  asymptotic upper and lower exponential bounds on the values that $\mu_\epsilon$ assigns to measurable subsets of $\mathbb{R}^n$. 

\begin{defn}
A rate function $I$ is a lower semicontinuous map 
 $I \, : \, \mathbb{R}^n \rightarrow [0,\infty]$, and it is a good rate function if for all $\alpha \in [0,\infty)$, the level set $\Psi_{I}(\alpha):= \{x \, : \, I(x) \leq \alpha\}$ is compact.
\end{defn}

For any 
set $B\subseteq\R^n$, we denote by  $B^{\circ}$ the interior of $B$. 

\begin{defn}\label{largdev}
A family of probability measures $\{\mu_\epsilon\}$ satisfies the large deviation principle with a rate function $I$ if, for all $B \in \mathcal{B}$,
\begin{equation}\label{lar}
-\inf_{x \in B^{\circ}} I(x) \leq \liminf_{\epsilon \rightarrow 0} \epsilon \log \mu_\epsilon (B) \leq \limsup_{\epsilon \rightarrow 0} \epsilon \log\mu_\epsilon(B) \leq - \inf_{x \in \bar{B}} I(x).
\end{equation}
\end{defn}

The right-and left-hand sides of  \eqref{lar} are referred to as the upper and lower bounds, respectively.

\begin{defn}
A family of probability measures $\{\mu_\epsilon\} $ on $\mathbb{R}^n$ is exponentially tight if for every $\alpha < \infty$, there exists a compact set $K_\alpha \subset \mathbb{R}^n$ such that
$$
\limsup_{\epsilon \rightarrow 0} \epsilon \log \mu_\epsilon(K_\alpha^c)< -\alpha.
$$
\end{defn}
Moreover, for each Borel measurable function $h: \mathbb{R}^n \rightarrow \mathbb{R}$, define
$$
\Lambda^{\epsilon}_h:=  \epsilon \log \int_{\mathbb{R}^n} e^{\frac{h(x)}{\epsilon}} \mu_\epsilon(dx).
$$
and  
\begin{equation}\label{limit}
\lim_{\epsilon \rightarrow 0} \epsilon \log \int_{\mathbb{R}^n} e^{\frac{h(x)}{\epsilon}} \mu_\epsilon(dx)=\Lambda_h
\end{equation}
provided the limit exists. Then, the so-called Bryc's inverse Varadhan Lemma permits to derive the large deviation principle as a consequence of exponential tightness of the measures $\mu_\epsilon$ and the existence of the limits  \eqref{limit} for every $h \in BC(\mathbb{R}^n)$. The statement is the following.

\begin{lem}\label{lem:bryc}
Suppose that the family $\{\mu_\epsilon\}$ is exponentially tight and that the limit in \eqref{limit} exists for every $h \in BC(\mathbb{R}^n)$. Then $\{\mu_\epsilon\}$ satisfies the LDP with the good rate function
$$
I(x)=\sup_{h \in BC(\mathbb{R}^n)} \{h(x)-\Lambda_h\}.
$$
Furthermore, for every $h \in BC(\mathbb{R}^n)$,
$$
\Lambda_h=\sup_{x \in \mathbb{R}^n} \{h(x)-I(x)\}.
$$
\end{lem}

Finally we recall the optional sampling theorem. For further details see \cite{W}. 
\begin{thm}\label{thm:opt}
Let $M = \{M_t\}_{t\geq 0}$ be a submartingale 
right-continuos and let $\tau$ be a stopping time, such that one of the following conditions is satisfied
\begin{itemize}
\item
$\tau$ is a.s. bounded, i.e. there exists $T\in (0,\infty)$ such that $\tau \leq T$ a.s.;
\item
$\tau$ is a.s. finite and $M_{\tau \wedge t} \leq Y$ for all $t \geq 0$, where $Y$ is an integrable variable (in particular $|M_{\tau \wedge n}| \leq K$ for a constant $K \in [0, \infty)$) 
\end{itemize}
Then the variable $M_\tau$ is integrable and
\begin{equation}\label{optio}
E(M_\tau) \geq E(M_0).
\end{equation}
If, instead,  $M$ is a supermartingale, then
$$
E(M_\tau) \leq E(M_0).
$$
\end{thm}

\end{document}